\documentclass[10pt]{article}
\usepackage{amsfonts}
\usepackage{amsmath}
\numberwithin{equation}{section}
\newtheorem{thm}{Theorem}[section]
\newtheorem{definition}{Definition}[section]

\newtheorem{propo}{Proposition}[section]
\newtheorem{remark}{Remark}[section]

\newtheorem{lemma}{Lemma}[section]
\newtheorem{coro}{Corollary}[section]
\newenvironment{proof}[1][Proof]{\noindent\textbf{#1.} }{\ \rule{0.5em}{0.5em}}
\newcounter{num}
\def\N{\Bbb N}
\def\Z{\Bbb Z}

\def\R{\Bbb R}
\def\C{\Bbb C}

\def\ds{\displaystyle}
\begin{document}
\title{\bf A Fractional Power for
Dunkl Transforms in $ L^{2}(\R^{N}, \ \omega_{k}(x)dx)$.}
\scrollmode
\author{Sami Ghazouani
\thanks{Institut Pr\'eparatoire aux Etudes d'Ing\'enieur de Bizerte,
Universit\'e de Carthage, 7021 Jarzouna, Tunisie. E-mail:
Ghazouanis@yahoo.fr.} \ and Fethi Bouzaffour
\thanks{Department of mathematics, College of Sciences, King Saud University, P. O Box 2455 Riyadh 11451, Saudi Arabia. E-mail : fbouzaffour@ksu.edu.sa.}}
\date{ }
\maketitle \vspace{0mm}

\begin{abstract}
A new fractional version of the Dunkl transform for real order $\alpha$ is obtained. An integral representation, a Bochner type identity and a Master formula for this transform are derived.
\end{abstract}
$\overline{\hspace{1.5cm}}$ \par Keywords: Fractional Fourier
transform, fractional Hankel transform,  Fractional Dunkl transform, Generalized Hermite polynomials and functions, semigroups of operators.
\section{Intorduction}
In recent years, there has been considerable interest in fractional versions of classical integral transform. For example the Fractional Fourier Transform (FFT), which may be considered as a fractional generalization of the classical Fourier Transform, becomes a remarkably powerful tool in signal processing, optics and quantum mechanics \cite{Almeida2}, \cite{Almeida1}, \cite{Alieva1} \cite{Alieva} and \cite{Namias}. The idea of Fractional powers of the Fourier Transform operator appeared firstly in the Wiener's works \cite{Winer}. The development of a wide-ranging modern theory, including operational formula, stems from a paper by Namias \cite{Namias}. As Namias's innovative ideas and results were developed in a formal manner, they were later revisited by McBride and Kerr \cite{McBride}, where a mathematically rigorous account is presented for the FFT on the space $\mathcal{S}(\R^{N})$ of rapid descent functions. More recently, Zayed \cite{Zayed} used the same approach as Namias to produce more general
method in order to defining fractional versions of a wider transforms class such that the fractional Hankel transform \cite{Kerr1}, \cite{Kerr}, \cite{Namias1}, the fractional integration and differentiation
operators.

This paper deals with the constraction of a fractional power of the Dunkl transform called the fractional Dunkl transform (FDT), using the Zayed's approach \cite{Zayed} and the multivariable generalized Hermite function introduced by R\"{o}sler \cite{rosler}. The resulting family of operators $\left\{D_{k}^{\alpha}\right\}_{\alpha \in \R}$ was proved to be a $\mathcal{C}_{0}$-group of unitary operators on $L^{2}\left(\R^{N},  \omega_{k}(x)\ dx\right),$ with infinitesimal generator $T.$ The spectral properties of $T$ is studied using the semigroup techniques. The FDT given in this paper has an integral representation which used with the analogue of the Funk-Hecke formula for $k$-spherical harmonics \cite{yuanxu} to derive a Bochner type identity for the FDT. The Master formula of the FDT is proved and founded to be generalizing the one given by  R\"{o}sler \cite{rosler} in Proposition 3.10. This Master formula is used to develop a new proof of the statement $(2)$ of the 3.4 R\"{o}sler's theorem \cite{rosler}.

The contents of the present paper are as follows. In section $2$, some basic definitions and results about harmonic analysis associated with Dunkl operators are collected. In section $3,$ the fractional Dunkl transform definition is given and then some elementary properties of this transformation are listed. The spectral properties of $T$ is studied. In section $4,$ the integral representation of the fractional Dunkl transform as well as the Bochner type identity and Master formula are given. In section $5,$ we find a subspace of $L^{2}(\R^{N}, \ \omega_{k}(x)dx)$ in which we define $T$ explicitly.
\section{Background: Dunkl operator }
In this section, we recall some notations and results on Dunkl
operators, Dunkl transform, and generalized Hermite functions   (see, \cite{Dunkl1,
Dunkl2, deJeu, opdam, rosler1}). In $\mathbb{R}^N$, we consider the
standard inner product
$$\langle x,\,y\rangle=\sum_{k=1}^Nx_ky_k$$
and the norm $|x|=\sqrt{\langle x,\,x\rangle}.$\\ For $u \in
\R^{N}\backslash \{0\}$, let $\sigma_{u}$ be the reflection on the
hyperplane $(\R u)^{\perp}$ orthogonal to $u$. It is given by
\begin{eqnarray}
\sigma_{u}(x)=x-2\ds\frac{\langle u,\,x\rangle
}{|u|^{2}}u. \label{x1}
\end{eqnarray}
A root system is a finite set $R$ of nonzero vectors
in $\mathbb{R}^N$ such that  $$R\cap \mathbb{R}u=\{\pm u\}\,\,\,\,
\text{and} \,\,\,\sigma_{u}R=R,$$ for all $u \in R.$

For a given root system $R$ normalized such that
$\langle u,u\rangle =2$ for all $u \in R$, we denote by $W$
the subgroup of $ O(N,\R)$ generated by the reflections $\{\sigma_{u}: \ u \in R\}$ and for $v \in
\R^{N}\backslash \underset{u\in R}{\ds\cup}<u>^{\perp},$ we fix the positive
subsystem $R_{+}=\{u \in R / \langle u,\,v \rangle >0 \}.$

A function $k:R\rightarrow \C$ on the root system $R$ is called a multiplicity
function, if for every $u \in R$ and $\sigma \in W, \  k(\sigma
u)=k(u).$

The Dunkl operators $T_{\xi}:=T_{\xi}(k),\xi \in \R^{N},$ associated with the reflection group $W$ and the multiplicity function $k$ are given by
\begin{eqnarray}
T_{\xi}f(x):=\partial_{\xi}f(x)+\ds\sum_{\eta \in R_{+}}k(\eta)<\eta,\xi>\ds\frac{f(x)-f(\sigma_{\eta}x)}{\langle \eta , x \rangle },\label{x2}
\end{eqnarray}
where $\partial_{\xi}$ denotes the derivative in the direction of
$\xi.$ Thanks to the $W$-invariance of the function $k$, this definition is independent of the choice of
the positive subsystem $R_{+}.$ For the $i$-th standard basis vector $\xi=e_{i} \in \R^{N}$, we use the abbreviation $T_{i}=T_{e_{i}}.$ The most striking property of Dunkl operators $T_{\xi}$, which is the foundation for rich analytic structures with them, is the following \cite{Dunkl1}:
\begin{thm}
For fixed $k$, the associated $T_{\xi}=T_{\xi}(k)$, $\xi \in \R^{N}$ commute.
\end{thm}
For fixed non-negative multiplicity functions $k$, let $\omega_{k}$ be
the weight function on $\R^{N}$ defined by
\begin{eqnarray*}
\omega_{k}(x)=\ds\prod_{\eta \in R_{+}}|\langle \eta , x\rangle|^{2k(\eta)}.
\end{eqnarray*}
It is $W$-invariant and homogeneous of degree $2\gamma,$ with the index
\begin{eqnarray*}
\gamma=\gamma(k)=\ds\sum_{\eta \in R_{+}}k(\eta).
\end{eqnarray*}

\textbf{Notation}:
We denote by $\Z_{+}$ the set of non-negative integers. For a multi-index $\nu=(\nu_{1},\dots,\nu_{N})\in \Z_{+}^{N},$ we write $|\nu|=\nu_{1}+\dots+\nu_{N}.$ The $\C$-algebra of polynomial functions on $\R^{N}$ is denoted by $\mathcal{P}=\C[\R^{N}].$
It has a natural grading
$$\mathcal{P}=\ds\bigoplus_{n\geq 0}\mathcal{P}_{n},$$
where $\mathcal{P}_{n}$ is the subspace of homogenous polynomials of (total) degree $n$. $\mathcal{S}(\R^{N})$ is the Schwartz space of rapidly decreasing functions on $\R^{N}.$
Finally, $L^{p}\left(\R^{N}, \ \omega_{k}(x)\ dx\right)$ is the space of measurable functions on $\R^{N}$ such that
\begin{eqnarray*}
\| f \|_{p}&=&\left(\ds\int_{\R^{N}}|f(y)|^{p}\omega_{k}(y)\ dy\right)^{\frac{1}{p}}< +\infty, \quad \mbox{if} \ 1\leq p<+\infty.
\end{eqnarray*}

The generalized Laplacian associated with $W$ and $k$, is defined by $\Delta_{k}:=\ds\sum_{i=1}^{N}T_{i}^{2}$. It is homogeneous of degree $-2$ on $\mathcal{P}$ and given explicitely by
\begin{equation}\Delta_{k}f(x)=\Delta f(x)+2\ds\sum_{\alpha \in R} k(\alpha)\left[\ds\frac{\left\langle \nabla f(x),\alpha\right\rangle}{\langle \alpha , x\rangle }-\ds\frac{f(x)-f(\sigma_{\alpha}x)}{\langle \alpha,x\rangle ^{2}}\right].\label{laplace}
\end{equation}
(Here $\Delta$ and $\nabla$ denote the usual Laplacian and gradient respectively).\\ \mbox{}\\

For $y\in \C^{N},$ the system
\begin{eqnarray}
\left\{\begin{array}{ll}
  T_{\xi}f(x,y)=\langle \xi,y\rangle f(x,y), & \mbox{for all} \ \xi \in\R^{N} , \\
f(0,y)=1, & \
\end{array}\right. \label{x3}
\end{eqnarray}
admits a unique analytic solution $K( \ . \ ,y)$ on $\R^{N}, $ which
is called the Dunkl kernel. Moreover, $K(.,.)$ extends to holomorphic function on $\C^{N}\times\C^{N}$
and possesses the following properties:\\
$1)$ For every $x,y\in \C^{N},$ we have
\begin{eqnarray}
K(x,y)=K(y,x) \quad \mbox{and } \quad K(\lambda x,y)=K(x,\lambda y) \ \mbox{for all} \ \lambda \in \C.\label{xx}
\end{eqnarray}
$2)$ For every $x,y\in \C^{N}$ and $g\in W,$ we have
\begin{eqnarray}
\overline{K( x,y)}=K( \overline{x},\overline{y})\quad \mbox{and } \quad K(gx,gy)=K( x,y).
\end{eqnarray}
$3)$ For all $\nu=(\nu_{1},\dots,\nu_{N}) \in \Z_{+}^{N}, \ x\in \R^{N} $ and $y\in \C^{N},$ we have
\begin{eqnarray}
\left|D_{y}^{\nu}K(x,y)\right|\leq |x|^{|\nu|}e^{|x||\Re y|}, \label{uu}
\end{eqnarray}
and for all $x,y\in \R^{N}:$
\begin{eqnarray}
\left|K(ix,y)\right|\leq 1, \label{v}
\end{eqnarray}
with $D_{y}^{\nu}=\ds\frac{\partial^{\nu}}{\partial y_{1}^{\nu_{1}}\dots\partial y_{N}^{\nu_{N}}}.$

The Dunkl transform associated with $R$ and $k\geq 0$ is defined
on $L^{1}\left(\R^{N}, \ \omega_{k}(x)\ dx\right)$ by
\begin{eqnarray*}
D_{k}f(x)=\ds\frac{c_{k}}{2^{\gamma+N/2}}\ds\int_{\R^{N}}f(y)K(-ix,y) w_{k}(y)dy,\label{}
\end{eqnarray*}
where the Mehta-type constant $c_{k}$ is given by
\begin{eqnarray*}
c_{k}=\left(\ds\int_{\R^{N}}e^{-|x|^{2}}w_{k}(x) dx \right)^{-1}. \label{}
\end{eqnarray*}
Some of the properties of the Dunkl transform are collected below
\cite{deJeu}
\begin{thm} \mbox{}\\
$a)$ (\textbf{Riemann-Lebesgue lemma}) For all $f \in L^{1}\left(\R^{N}, \ \omega_{k}(x)\ dx\right),$ the Dunkl transform $D_{k}f$ belongs to $C_{0}(\R^{N}).$\\
$b)$ (\textbf{$L^{1}$-inversion}) For all $f \in L^{1}\left(\R^{N}, \ \omega_{k}(x)\ dx\right)$ with $D_{k}f \in L^{1}\left(\R^{N}, \ \omega_{k}(x)\ dx\right)$,
\begin{eqnarray}
D_{k}^{2}f=\check{f}, \  a.e, \ \mbox{where} \ \check{f}(x)=f(-x).  \label{w3}
\end{eqnarray}
$c)$ The Dunkl transform $f\rightarrow D_{k}f$ is an automorphism of $\mathcal{S}(\R^{N}).$
\\
$d)$ (\textbf{Plancherel Theorem})\\
$i)$ If $f\in L^{1}(\R^{N},  \omega_{k}(x)\ dx)\cap L^{2}(\R^{N},  \omega_{k}(x)\ dx),$ then $D_{k}f \in L^{2}(\R^{N},  \omega_{k}(x)\ dx)$ and $\| D_{k}f \|_{2}=\| f \|_{2}.$\\
$ii)$ The Dunkl transform has a unique extension to an isometric isomorphism of $L^{2}\left(\R^{N},  \omega_{k}(x)\ dx\right).$ The extension is also denoted by $f\rightarrow D_{k}f.$
\label{tt1}
\end{thm}
In \cite{rosler}, R\"{o}sler constructed systems of naturally
associated multivariable generalized Hermite polynomials $\{H_{\nu},\nu \in \Z_{+}^{N}\}$ and Hermite functions $\{h_{\nu}, \nu \in \Z_{+}^{N}\},$
\begin{eqnarray}
H_{\nu}(x)&=&2^{|\nu|}e^{-\Delta_{k}/4}\varphi_{\nu}(x), \\ h_{\nu}(x)&=& \ds\frac{\sqrt{c_{k}}}{2^{|\nu|/2}} \ e^{-|x|^{2}/2}H_{\nu}(x), \label{x}
\end{eqnarray}
such that $\varphi_{\nu} \in \mathcal{P}_{|\nu|}$ and the coefficients of $\varphi_{\nu}$'s are real. He proved that the generalized Hermite functions $\{ h_{\nu}, \ \nu \in \Z_{+}^{N}\}$ form an an orthonormal basis of eigenfunctions for the Dunkl operator on $L^{2}(\R^{N},\ \omega_{k}(x) dx)$ with
\begin{eqnarray}
D_{k}h_{\nu}= (-i)^{|\nu|}h_{\nu}. \label{nn}
\end{eqnarray}
\section{Fractional Dunkl transforms in $L^{2}(\R^{N}, \ \omega_{k}(x)dx)$.}
In the remainder of this paper, we denote by $R$ a root system in $\mathbb{R}^N,$ $R_+$ a fixed positive
subsystem of $R$ and $k$ a nonnegative multiplicity function defined on
$R.$
\subsection{Definition and properties.}
Consider the Hilbert space $L^{2}(\R^{N}, \ \omega_{k}(x)dx)$ equipped with the orthonormal basis $\{h_{\nu}, \nu \in \Z_{+}^{N}\}.$\\
Let $f \in L^{2}(\R^{N}, \ \omega_{k}(x)dx),$ then we can expand $f$ in the orthonormal basis $\{h_{\nu}, \nu \in \Z_{+}^{N}\}$ as follows:
$$f=\ds\sum_{\nu \in \Z_{+}^{N}}\langle f,h_{\nu} \rangle  \ h_{\nu}=\ds\sum_{\nu \in \Z_{+}^{N}}\hat{f}_{\nu}  \ h_{\nu},$$
where
$$\hat{f}_{\nu}=\langle f, h_{\nu}\rangle=\ds\int_{\R^{N}}f(x)\overline{h_{\nu}(x)}w_{k}(x)dx.$$
From Theorem \ref{tt1}, $d), ii)$ and (\ref{nn}), the identity
\begin{eqnarray}
D_{k}f&=&\ds\sum_{\nu \in \Z_{+}^{N}} (-i)^{|\nu|} \langle f,h_{\nu} \rangle  \ h_{\nu}\nonumber\\
&=& \ds\sum_{\nu \in \Z_{+}^{N}} e^{-i\frac{\pi}{2}|\nu|}\langle f,h_{\nu} \rangle  \ h_{\nu},\label{oo}
\end{eqnarray}
holds for every $f \in L^{2}(\R^{N}, \ \omega_{k}(x)dx),$ which allows us to give the following definition:
\begin{definition}
Let $\alpha \in \R,$ we define the fractional Dunkl transform $D_{k}^{\alpha}$ by
\begin{eqnarray}
D_{k}^{\alpha}f&=& \ds\sum_{\nu \in \Z_{+}^{N}} \ e^{i|\nu|\alpha}  \  \langle f,h_{\nu} \rangle  \ h_{\nu} . \label{a}
\end{eqnarray}
\end{definition}
\begin{remark}The limit in the second member of (\ref{a}) exists. In fact, let us denote by $\mathcal{P}_{0}(\Z_{+}^{N})$ the set of finite subsets of $\Z_{+}^{N}.$ For every $J\subset\mathcal{P}_{0}(\Z_{+}^{N}),$ we have
\begin{eqnarray}
\left\|\sum_{\nu \in J } \ e^{i|\nu|\alpha}  \  \langle f,h_{\nu} \rangle  \ h_{\nu}\right\|_{2}^{2}&=&\sum_{\nu \in J }  |\langle f,h_{\nu} \rangle|^{2}. \label{c}
\end{eqnarray}
The convergence in $L^{2}(\R^{N}, \ \omega_{k}(x)dx)$ of the series
$$\ds\sum_{\nu \in \Z_{+}^{N}} \ e^{i|\nu|\alpha}  \  \langle f,h_{\nu} \rangle  \ h_{\nu},$$
is obtained by means of the Cauchy convergence test and (\ref{c}).
\end{remark}
We summarize the elementary properties of $D_{k}^{\alpha}$ in the next  Proposition.
\begin{propo} Let $\alpha, \ \beta \in \R .$ The fractional Dunkl transform $D_{k}^{\alpha}$ satisfies the following properties: \\
$1)$ \ $D_{k}^{0}=I$, which is the identity operator, \\
$2)$ \ $D_{k}^{\pi/2}=D_{k}$, \\
$3)$ \ $D_{k}^{\alpha}\circ D_{k}^{\beta}=D_{k}^{\alpha+\beta}$, \\
$4)$ \ $D_{k}^{\alpha+2\pi}=D_{k}^{\alpha}$, \\
$5)$ \ $D_{k}^{\pi}=\hat{I}$, where $\hat{I}f(x)=f(-x)$,\\
$6)$ \ For all $f$ and $g \in L^{2}(\R^{N}, \ \omega_{k}(x)dx), \
\langle D_{k}^{\alpha}f,g\rangle=\langle f,D_{k}^{-\alpha}g\rangle.$
\label{pr02}
\end{propo}
\proof \\
$1)$ and $4)$ follow immediately from (\ref{a}).\\
$2)$ Follows immediately from (\ref{oo}).\\
$3)$ From (\ref{a}), we have
\begin{eqnarray*}
D_{k}^{\alpha}(D_{k}^{\beta}f)&=&\ds\sum_{\nu \in \Z_{+}^{N}} \ e^{i|\nu|\alpha}  \  \langle D_{k}^{\beta}f,h_{\nu} \rangle  \ h_{\nu}\\ &=& \ds\sum_{\nu \in \Z_{+}^{N}}e^{i|\nu|(\alpha+\beta)} \langle f,h_{\nu} \rangle  \ h_{\nu}=D_{k}^{\alpha+\beta}f.
\end{eqnarray*}
$5)$ From (Lemma 3.11,
\cite{rosler}), we conclude that
$$ e^{-|x|^{2}/2} \ h_{\nu}(x)=\sqrt{c_{k}2^{|\nu|}}\ds\int_{\R^{N}}K(x,-2iy)\varphi_{\nu}(y)d_{\mu}(y),$$
where $\varphi_{\nu}$ is a homogenous polynomials of order $|\nu|$.
Thus
 \begin{eqnarray}
h_{\nu}(-x)&=&e^{|x|^{2}/2}\sqrt{c_{k}2^{|\nu|}}\ds\int_{\R^{N}}K(-x,-2iy)\varphi_{\nu}(y)d_{\mu}(y),\nonumber\\&=&
e^{|x|^{2}/2}\sqrt{c_{k}2^{|\nu|}}\ds\int_{\R^{N}}K(x,2iy)\varphi_{\nu}(y)d_{\mu}(y),\nonumber\\&=&(-1)^{|\nu|}
e^{|x|^{2}/2}\sqrt{c_{k}2^{|\nu|}}\ds\int_{\R^{N}}K(x,-2iy)\varphi_{\nu}(y)d_{\mu}(y),\nonumber\\&=&(-1)^{|\nu|}h_{\nu}(x).\nonumber
\end{eqnarray}
Hence,
\begin{eqnarray*}
D_{k}^{-\pi}f&=&\sum_{\nu \in Z_{+}^{N} } \ e^{-i|\nu|\pi}  \
\langle f,h_{\nu} \rangle  \ h_{\nu}=\sum_{\nu \in
Z_{+}^{N} } \ (-1)^{|\nu|}  \  \langle f,h_{\nu} \rangle  \
h_{\nu},\nonumber\\&=&\hat{I}f.\nonumber
\end{eqnarray*}
$6)$ Let $f$ and $g \in L^{2}(\R^{N}, \ \omega_{k}(x)dx).$ It is easy to check that
\begin{eqnarray*}
\langle D_{k}^{\alpha}f,g\rangle&=&\ds\sum_{\nu \in \Z_{+}^{N}} \ e^{i|\nu|\alpha}  \  \langle f,h_{\nu} \rangle  \ \overline{\langle g,h_{\nu} \rangle}=\ds\sum_{\nu \in \Z_{+}^{N}} \langle f,h_{\nu} \rangle  \ \overline{e^{-i|\nu|\alpha}\langle g,h_{\nu} \rangle}\\&=&
\langle f,D_{k}^{-\alpha}g\rangle.
\end{eqnarray*}
\begin{thm} The family of operators $\{D_{k}^{\alpha}\}_{\alpha \in \R}$ is a $\mathcal{C}_{0}$-group of unitary operators on $L^{2}(\R^{N},\ \omega_{k}(x)\ dx).$
\end{thm}
\proof \\ From Proposition \ref{pr02}, we deduce that
the family $\{D_{k}^{\alpha}\}_{\alpha \in \R}$
satisfies the algebraic properties of a group:
\begin{eqnarray*}
 D_{k}^{0}=I, \quad D_{k}^{\alpha}\circ D_{k}^{\beta}=D_{k}^{\alpha+\beta}=D_{k}^{\beta}\circ D_{k}^{\alpha} \ ; \quad \alpha, \ \beta \in \R.
\end{eqnarray*}
For the strong continuity, assume that
$f \in L^{2}(\R^{N},\ \omega_{k}(x)\ dx).$ Then
\begin{eqnarray*}
\left\|D_{k}^{\alpha}f-f\right\|_{2}^{2}=\ds\sum_{\nu \in \Z_{+}^{N}}|e^{i|\nu|\alpha}-1|^{2}\left|\langle f, h_{\nu}\rangle \right|^{2}.
\end{eqnarray*}
For each $\nu \in \Z_{+}^{N},$ we have
\begin{eqnarray*}
\ds\lim_{\alpha\rightarrow 0}|e^{i|\nu|\alpha}-1|^{2}\left|\langle f, h_{\nu}\rangle \right|^{2}&=&0,\\
|e^{i|\nu|\alpha}-1|^{2}\left|\langle f, h_{\nu}\rangle \right|^{2}&\leq& 4\left|\langle f, h_{\nu}\rangle \right|^{2}.
\end{eqnarray*}
Since
\begin{eqnarray*}
\ds\sum_{\nu \in \Z_{+}^{N}}\left|\langle f, h_{\nu}\rangle \right|^{2}=\|f\|_{2}^{2}<\infty,
\end{eqnarray*}
then we can interchange limits and sum to get:
\begin{eqnarray*}
\ds\lim_{\alpha\rightarrow 0}\left\|D_{k}^{\alpha}f-f\right\|_{2}^{2}= 0.
\end{eqnarray*}
Hence $\{D_{k}^{\alpha}\}_{\alpha \in \R}$ is
a strongly continuous group of operators on $L^{2}(\R^{N},\ \omega_{k}(x)\ dx).$ In addition, by Proposition \ref{pr02}, we have for all $f, \ g \in L^{2}(\R^{N},\ \omega_{k}(x)\ dx), $
\begin{eqnarray*}
\langle D_{k}^{\alpha}f, \ g\rangle =\langle f, \ D_{k}^{-\alpha}g\rangle,
\end{eqnarray*}
and therefore $(D_{k}^{\alpha})^{*}=D_{k}^{-\alpha}=(D_{k}^{\alpha})^{-1},$ establishing that each $D_{k}^{\alpha}$ is unitary.\\
We have therefore shown that $\{D_{k}^{\alpha}\}_{\alpha \in \R}$ is a $\mathcal{C}_{0}$-group of unitary operators on $L^{2}(\R^{N}, \ \omega_{k}(x)dx).$ \\
The infinitesimal generator $T$ of $\{D_{k}^{\alpha}\}_{\alpha \in \R}$ is defined by
$$
\begin{tabular}{rrl}
$T: L^{2}(\R^{N}, \ \omega_{k}(x)dx)\supseteq D(T)$& $\longrightarrow$ & $L^{2}(\R^{N}, \ \omega_{k}(x)dx),$\\
    $f$ & $\longmapsto$ & $Tf$ \\
\end{tabular}
$$
where
\begin{eqnarray*}
D(T)&=&\left\{f\in L^{2}(\R^{N}, \ \omega_{k}(x)dx): \ \ds\lim_{\alpha\rightarrow 0}(1/\alpha)[D_{k}^{\alpha}f-f]\in L^{2}(\R^{N}, \ \omega_{k}(x)dx)  \right\},\\
  Tf &=&\ds\lim_{\alpha\rightarrow 0}(1/\alpha)[D_{k}^{\alpha}f-f], \quad f \in D(T).
\end{eqnarray*}
\subsection{Spectral properties of the operator $T$.}
We denote by $\mathcal{B}(L^{2}(\R^{N}, \ \omega_{k}(x)dx)),$ the set of all linear bounded operator in $L^{2}(\R^{N}, \ \omega_{k}(x)dx).$
The resolvent set of $T$ is the set $\rho(T)$ consisting of all scalars $\lambda$ for which the linear operator $\lambda I-T$ is a 1-1 mapping from its domain $D(\lambda I-T)=D(T)$ on to the Hilbert space $L^{2}(\R^{N}, \ \omega_{k}(x)dx)$ with $(\lambda I-T)^{-1}\in \mathcal{B}(L^{2}(\R^{N}, \ \omega_{k}(x)dx)).$ The spectrum of $T$ is the set $\sigma(T)$ that is the complement of $\rho(T)$ in $\C.$ The function $R(\lambda,T)=(\lambda I-T)^{-1}$
from $\rho(T)$ into $\mathcal{B}(L^{2}(\R^{N}, \ \omega_{k}(x)dx))$ is the resolvent of $T.$

As $T$ is the generator of the $\mathcal{C}_{0}$-group $\{D_{k}^{\alpha}\}_{\alpha \in \R},$ some elementary properties of $T$ and $D_{k}^{\alpha}$ are listed in the following proposition (see \cite{Engel}, \cite{Goldstein}).
\begin{propo} Let $\alpha \in \R.$ The following properties hold.\\
$i)$ If $f\in D(T),$ then $D_{k}^{\alpha}f \in D(T)$ and
\begin{eqnarray}
\frac{d}{d\alpha} D_{k}^{\alpha}f=D_{k}^{\alpha}Tf=TD_{k}^{\alpha}f.
\end{eqnarray}
$ii)$ For every $t\in \R$ and $f\in L^{2}(\R^{N}, \ \omega_{k}(x)dx),$ one has
$$\ds\int_{0}^{t}D_{k}^{\alpha}f \ d\alpha \in D(T).$$
$iii)$ For every $\alpha\in \R,$ one has
\begin{eqnarray}
D_{k}^{\alpha}f-f&=&T\ds\int_{0}^{\alpha}D_{k}^{s}f \ ds, \quad \mbox{if} \ f \in L^{2}(\R^{N}, \ \omega_{k}(x)dx)\\
&=& \ds\int_{0}^{\alpha}D_{k}^{s}Tf \ ds, \quad \mbox{if} \ f \in D(T).
\end{eqnarray}
\label{pro3}
\end{propo}
\begin{remark} If we apply the Proposition \ref{pro3}, $iii)$ to the rescaled semigroup
$$S(\alpha):=e^{-\lambda \alpha}D_{k}^{\alpha}, \quad \alpha\in \R$$
\end{remark}
whose generator is $B:=T-\lambda I $ with domain $D(B)=D(T),$ we obtain
for every $\lambda \in \C$ and $ \alpha\in \R,$
\begin{eqnarray}
-e^{-\lambda \alpha}D_{k}^{\alpha}f+f &=&(\lambda I-T)\ds\int_{0}^{\alpha}e^{-\lambda s}D_{k}^{s}f \ ds; \quad \ f\in L^{2}(\R^{N}, \ \omega_{k}(x)dx),\label{b}\\
&=&\ds\int_{0}^{\alpha}e^{-\lambda s}D_{k}^{s}(\lambda I-T)f \ ds; \quad \ f\in D(T).\label{d}
\end{eqnarray}
Now we are interesting with the eigenvalues of $T$ by giving an important formula relating the semigroup $\{D_{k}^{\alpha}\}_{\alpha \in \R},$ to the resolvent of its generator $T.$
\begin{propo} For the operator $T,$ the following properties hold:\\
$1)$ T is closed and densely defined.\\
$2)$ The operator $iT$ is self-adjoint.\\
$3)$ $\sigma(T)=\sigma_{p}(T)\subset i \Z,$ and for each $\lambda \in \C\backslash i\Z $ and for all $f \in L^{2}(\R^{N}, \ \omega_{k}(x)dx),$
\begin{eqnarray}
R(\lambda,T)f=(1-e^{-2\pi \lambda})^{-1}\ds\int_{0}^{2\pi}e^{-\lambda s}D_{k}^{s}f \ d s. \label{l}
\end{eqnarray}
Here $\sigma_{p}(T)=\left\{\lambda \in \C: \ \lambda I-T \ \mbox{is not injective} \right\}.$ \\
\end{propo}
\proof\\
$1)$ The fact that  T is closed and densely defined follows from the Hille-Yosida Theorem (see[\cite{Goldstein}, p. 15]).\\
$2)$ Since $\{D_{k}^{\alpha}\}_{\alpha \in \R}$ is unitary, it follows from Stone's Theorem [\cite{Goldstein}, p. 32] that $T$ is skew-adjoint ($T^{*}=-T$) and therefore $iT$ is self-adjoint.\\
$3)$ If we replace $\alpha$ by $2\pi$ in (\ref{b}) and (\ref{d}) and we use the fact that $\{D_{k}^{\alpha}\}_{\alpha \in \R}$ is a periodic $\mathcal{C}_{0}$-group with period $2\pi,$ we get
\begin{eqnarray}
(1-e^{-2\pi\lambda })f &=&(\lambda I-T)\ds\int_{0}^{2\pi}e^{-\lambda s}D_{k}^{s}f \ ds; \quad \ f\in L^{2}(\R^{N}, \ \omega_{k}(x)dx),\label{m}\\
&=&\ds\int_{0}^{2\pi}e^{-\lambda s}D_{k}^{s}(\lambda I-T)f \ ds; \quad \ f\in D(T).\label{n}
\end{eqnarray}
Let $\lambda \not\in i\Z.$ Then  $1-e^{-2\pi\lambda }\neq 0.$ By the use of (\ref{m}) and (\ref{n}), $\lambda I -T$ is invertible ($\lambda \in \rho(T)$) and
$$(\lambda I-T)^{-1}f=R(\lambda,T)f=(1-e^{-2\pi \lambda})^{-1}\ds\int_{0}^{2\pi}e^{-\lambda s}D_{k}^{s}f \ d s.$$
The previous Proposition indicates that every point in the spectrum of $T$ is an isolated point of the set $i\Z.$ Let $i n$ be an element of the spectrum of $T$ and
$$P_{n}=\ds\frac{1}{2i\pi}\ds\int_{\Gamma}R(\lambda, T)\ d\lambda,$$
the associated spectral projection, where $\Gamma$ is a Jordan path in the complement of $i\Z\backslash \{in\}$ and enclosing $in.$ The function $\lambda\longmapsto R(\lambda,T)$ can be expanded as a Laurent series
$$R(\lambda,T)=\ds\sum_{k=-\infty}^{+\infty}(\lambda-in)^{k}B_{k}$$
for $0 < |\lambda-in|<\delta$ and some sufficiently small $\delta>0.$ The coefficients $B_{k}$ of this series are bounded operators given by the formulas
$$B_{k}=\ds\frac{1}{2i\pi}\ds\int_{\Gamma}\ds\frac{R(\lambda, T)}{(\lambda-in)^{k+1}}\ d\lambda, \quad k \in \Z.$$
The coefficient $B_{-1}$ is exactly the spectral projection $P_{n}$ corresponding to the decomposition $\sigma(T)=\{in\}\cup \{i\Z\backslash \{in\}\}$ of the spectrum of $T.$ From (\ref{l}), one deduces the identitie
\begin{eqnarray}
P_{n}&=&B_{-1}=\ds\lim_{\lambda\rightarrow in}(\lambda-in)R(\lambda,T)\nonumber\\
&=&\frac{1}{2\pi}\ds\int_{0}^{2\pi}e^{-in s}D_{k}^{s} \ ds,\label{yy}
\end{eqnarray}
which allows as to interpret $P_{n}$ as the nth Fourier coefficient of the $2\pi$-periodic function $s\longmapsto D_{k}^{s}.$

In the following Proposition we gather some properties of the operator $P_{n}.$
\begin{propo} Let $n,m \in \Z$ such that $n\neq m$ and $f,g \in L^{2}(\R^{N}, \ \omega_{k}(x)dx).$ Then\\
$i)$ $TP_{n}=inP_{n},$\\
$ii)$ $D_{k}^{s}P_{n}=e^{ins}P_{n},$\\
$iii)$ $P_{n}P_{m}=0,$\\
$iv)$ $\left\langle P_{n}f,g\right>=\left\langle f,P_{n}g\right\rangle.$ In particular $\left\langle P_{n}f,P_{m}g\right>=0.$ \\
$v)$ The linear span $$\mbox{lin} \ \underset{n\in\Z}{\bigcup}P_{n}L^{2}(\R^{N}, \ \omega_{k}(x)dx)$$
is dense in $L^{2}(\R^{N}, \ \omega_{k}(x)dx).$\label{pro4}
\end{propo}
\proof \mbox{}\\
$i)$ It follows directly from (\ref{m}) applied to $\lambda=in.$\\
$ii)$ Applying $D_{k}^{s}$ to each member of (\ref{yy}) then according to Proposition \ref{pr02}, $3)$, we obtain
\begin{eqnarray*}
D_{k}^{s}P_{n}f&=&\frac{1}{2\pi}\int_{0}^{2\pi}e^{-int}D_{k}^{s}D_{k}^{t}f \ dt \\&=&\frac{1}{2\pi}\int_{0}^{2\pi}e^{-int}D_{k}^{s+t}f \ dt.
\end{eqnarray*}
The change of variables $u=s+t$ gives the desired result.\\
$iii)$ From (\ref{yy}) and $ii)$, we have
\begin{eqnarray*}
P_{n}P_{m}f&=&\frac{1}{2\pi}\int_{0}^{2\pi}e^{-ins}D_{k}^{s}(P_{m}f) \ ds \\&=&\left(\frac{1}{2\pi}\int_{0}^{2\pi}e^{i(m-n)s}  \ ds \right)\ P_{m}f\\&=& 0.
\end{eqnarray*}
$iii)$ Obvious.\\
$iv)$ Assume that the linear span
$$\mbox{lin} \ \underset{n\in\Z}{\bigcup}P_{n}L^{2}(\R^{N}, \ \omega_{k}(x)dx)$$
is not dense in $L^{2}(\R^{N}, \ \omega_{k}(x)dx).$  By the Hahn-Banach theorem there exists a nonzero linear functional
$$ \varphi:L^{2}(\R^{N}, \ \omega_{k}(x)dx)\longrightarrow \C$$
vanishing on each $P_{n}L^{2}(\R^{N}, \ \omega_{k}(x)dx), \ n\in \Z.$ By the Riesz representation theorem, there exists a unique vector $g \in L^{2}(\R^{N}, \ \omega_{k}(x)dx)\backslash\{0\}$ such that
$$\varphi(f)=\left\langle f,g\right\rangle \ \mbox{for all} \ f \in L^{2}(\R^{N}, \ \omega_{k}(x)dx).$$
Hence for all $n\in \Z $ and $f \in L^{2}(\R^{N}, \ \omega_{k}(x)dx),$
\begin{eqnarray*}
0=\left\langle P_{n}f,g\right\rangle &=&\left\langle \frac{1}{2\pi}\int_{0}^{2\pi}e^{-ins}D_{k}^{s}f \ ds, g\right\rangle \nonumber\\&=&\frac{1}{2\pi}\int_{0}^{2\pi}e^{-ins} \left\langle D_{k}^{s}f,g\right\rangle \ ds .
\end{eqnarray*}
For each $f\in L^{2}(\R^{N}, \ \omega_{k}(x)dx),$ the function $s\longmapsto \left\langle D_{k}^{s}f,g\right\rangle$ has all its Fourier coefficients equal to zero, then it vanishes. This cannot be true, since if we take $f=g$ and $s=0,$
$$\left\langle D_{k}^{0}g,g\right\rangle=\|g\|_{2}^{2}>0.$$
\begin{propo} Let $f\in D(T).$ Then
\begin{eqnarray}
f&=&\ds\sum_{n=-\infty}^{+\infty}P_{n}f,\label{bb}
\end{eqnarray}
and therefore, if $f\in D(T^{2})$
\begin{eqnarray}
Tf&=&\ds\sum_{n=-\infty}^{+\infty}inP_{n}f.\label{cc}
\end{eqnarray}\label{pro7}
\end{propo}
\proof
We are going to show that the series $\ds\sum_{n\in \Z}P_{n}f$ is summable for all $f\in D(T).$
For this, let $f\in D(T)$ and put $g=Tf.$ The commutativity of $T$ and $P_{n}$ together with Proposition \ref{pro4} gives:
 $$P_{n}g=P_{n}Tf=TP_{n}f=inP_{n}f.$$
By the Cauchy-Schwartz inequality, it follows that
 \begin{eqnarray*}
\left|\ds\sum_{n\in H}\left\langle P_{n}f,h\right\rangle\right|&=&\left|\ds\sum_{n\in H}(in)^{-1}\left\langle P_{n}g,h\right\rangle\right|\\&\leq& \left(\sum_{n\in H}n^{-2}\right)^{1/2}\left(\sum_{n\in H}\left|\left\langle P_{n}g,h\right\rangle\right|^{2}\right)^{1/2},
\end{eqnarray*}
where $h\in L^{2}(\R^{N}, \ \omega_{k}(x)dx) $ and $H$ be a finite subset of $\Z\backslash\{0\}.$
The function $s\longmapsto \left\langle D_{k}^{s}g,h\right\rangle$ belongs to $L^{2}([0,2\pi]),$ then we obtain from Bessel's inequality
\begin{eqnarray*}
\sum_{n\in H}\left|\left\langle P_{n}g,h\right\rangle\right|^{2}&\leq& \frac{1}{2\pi}\ds\int_{0}^{2\pi}|\langle D_{k}^{s}g,h\rangle|^{2} \ ds\\&\leq&\frac{ \|h\|_{2}^{2}}{2\pi}\ds\int_{0}^{2\pi} \|D_{k}^{s}g\|_{2}^{2} \ ds=\|h\|_{2}^{2} \ \|g\|_{2}^{2}.
\end{eqnarray*}
Therefore, for any $h\in L^{2}(\R^{N}, \ \omega_{k}(x)dx),$
\begin{eqnarray*}
\left|\left\langle \ds\sum_{n\in H}P_{n}f,h\right\rangle \right|=\left|\ds\sum_{n\in H}\left\langle P_{n}f,h\right\rangle\right|&\leq&\|h\|_{2} \ \|g\|_{2}\ \left(\sum_{n\in H}n^{-2}\right)^{1/2}.
\end{eqnarray*}
Taking supremum over $ h\in L^{2}(\R^{N}, \ \omega_{k}(x)dx)$ with $\| h\|_{2}\leq 1,$ we get
\begin{eqnarray*}
\left\| \ds\sum_{n\in H}P_{n}f\right\|_{2}&\leq& \|g\|_{2}\ \left(\sum_{n\in H}n^{-2}\right)^{1/2},
\end{eqnarray*}
which means that the series $\ds\sum_{n\in \Z}P_{n}f$ converges converges in $L^{2}(\R^{N}, \ \omega_{k}(x)dx).$\\ Set $$f_{1}=\ds\sum_{n=-\infty}^{+\infty}P_{n}f$$
and let $g\in L^{2}(\R^{N}, \ \omega_{k}(x)dx).$ As the Fourier coefficients of the continuous, $2\pi$-periodic functions
$$s\longmapsto \left\langle D_{k}^{s}f_{1},g\right\rangle \quad \mbox{and} \quad s\longmapsto \left\langle D_{k}^{s}f,g\right\rangle $$
coincide. Then, for all $s\in \R,$
$$\left\langle D_{k}^{s}f_{1},g\right\rangle=\left\langle D_{k}^{s}f,g\right\rangle.$$
In particular, for $s=0, \ \langle f_{1},g\rangle=\langle f,g\rangle$ and therefore $f_{1}=f.$ \\
Replacing $f$ in (\ref{bb}) by $Tf,$ then we get (\ref{cc}).\\
At the end of the section $4$, we will show that $P_{n}=0,$ for any negative integer $n\neq0.$
\section{The fractional Dunkl transform in $L^{1}(\R^{N}, \ \omega_{k}(x)dx)\cap L^{2}(\R^{N}, \ \omega_{k}(x)dx)$.}
In this section, we shall derive an integral representation for the fractional Dunkl transform $D_{k}^{\alpha}$ defined by (\ref{a}), for suitable function $f.$
\subsection{Integral representation.}
Following Zayed's approach \cite{Zayed}, we define the operator $D_{k,r}^{\alpha}$ as
\begin{eqnarray}
D_{k,r}^{\alpha}f:=\ds\sum_{\nu \in \Z_{+}^{N}}r^{|\nu|}e^{i|\nu|\alpha}\langle f, h_{\nu}\rangle h_{\nu},\label{}
\end{eqnarray}
where $0<r\leq1$ and so $D_{k}^{\alpha}=D_{k,1}^{\alpha}.$\\
In the next proposition, we collect some properties of $D_{k,r}^{\alpha}.$
\begin{propo} Let $\alpha \in \R $ and $r\in]0, 1].$ Then\\
$1)$ $D_{k,r}^{\alpha}$ is a bounded operator on $L^{2}(\R^{N}, \ \omega_{k}(x)dx)$ satisfying $\| D_{k,r}^{\alpha}f\|_{2}\leq \| f \|_{2}.$\\
$2)$ For all $f\in L^{2}(\R^{N}, \ \omega_{k}(x)dx), \ D_{k,r}^{\alpha}f\rightarrow D_{k}^{\alpha}f \ \mbox{in} \ L^{2}(\R^{N}, \ \omega_{k}(x)dx) \quad \mbox{as} \ r\rightarrow 1^{-}.
$
\end{propo}
\proof Let $f\in L^{2}(\R^{N}, \ \omega_{k}(x)dx).$\\
$1)$ According to Parseval's formula, we have
\begin{eqnarray*}
\| D_{k,r}^{\alpha}f\|_{2}^{2}&=&\ds\sum_{\nu \in \Z_{+}^{N}}r^{2|\nu|}\left|\langle f, h_{\nu}\rangle\right|^{2}\\
&\leq&  \ds\sum_{\nu \in \Z_{+}^{N}}\left|\langle f, h_{\nu}\rangle\right|^{2}=\| f\|_{2}^{2}.
\end{eqnarray*}
$2)$ It is easy to see that
$$D_{k,r}^{\alpha}f-D_{k}^{\alpha}f=\ds\sum_{\nu \in \Z_{+}^{N}}(r^{|\nu|}-1)e^{i|\nu|\alpha}\langle f, h_{\nu}\rangle h_{\nu}.$$
Then
$$\left\|D_{k,r}^{\alpha}f-D_{k}^{\alpha}f\right\|_{2}^{2}=\ds\sum_{\nu \in \Z_{+}^{N}}|r^{|\nu|}-1|^{2}\left|\langle f, h_{\nu}\rangle\right|^{2} .$$
By the dominated convergence theorem it follows that
$\ds\lim_{r\rightarrow 1^{-}}\left\|D_{k,r}^{\alpha}f-D_{k}^{\alpha}f\right\|_{2}^{2}=0.$
\begin{coro} For each fixed $f \in L^{2}(\R^{N}, \ \omega_{k}(x)dx),$ there exists $\{ r_{j}\}_{j=1}^{\infty},$ with $r_{j}\rightarrow 1^{-}$ as $j\rightarrow \infty,$ such that
\begin{eqnarray*}
D_{k}^{\alpha}f(x)=\ds\lim_{j\rightarrow\infty}D_{k,r_{j}}^{\alpha}f(x)
\end{eqnarray*}
for almost all $x \in \R^{N}.$\label{coro1}
\end{coro}
\proof \\ This is a consequence of a standard result that if a sequence $\{f_{n}\}$ converges in $L^{2}(\R^{N}, \ \omega_{k}(x)dx)$ to $f,$ then there exists a subsequence $\{f_{n_{k}}\}$ that converges pointwise almost everywhere to $f.$\\
The operator $D_{k,r}^{\alpha}$ defined above have the integral representation given in the next lemma.
\begin{lemma} For $f\in L^{2}(\R^{N}, \ \omega_{k}(x)dx)$ and $0<r<1,$ we have \begin{eqnarray}
D_{k,r}^{\alpha}f(x)=\ds\int_{\R^{N}}K_{\alpha}(r,x,y)f(y)\omega_{k}(y)\ dy,
\end{eqnarray}
where
\begin{eqnarray}
K_{\alpha}(r,x,y)=\ds\sum_{\nu\in \Z_{+}^{N}}r^{|\nu|}e^{i|\nu|\alpha}h_{\nu}(x)\overline{h_{\nu}(y)}. \label{e}
\end{eqnarray}
\end{lemma}
\proof Let $x\in \R^{N}$ and $H$ be a finite subset of $\Z_{+}^{N}$. Then
\begin{eqnarray}
\left\| \ds\sum_{\nu\in H}h_{\nu}(x)h_{\nu}(y)(re^{i\alpha})^{|\nu|} \right\|_{2}^{2}=\ds\sum_{\nu\in H}|h_{\nu}(x)|^{2}|r|^{2|\nu|} . \label{vv}
\end{eqnarray}
Since the series (see Theorem 3.12 in \cite{rosler})
$$\ds\sum_{\nu\in\Z_{+}^{N}}h_{\nu}(x)h_{\nu}(y)(re^{i\alpha})^{|\nu|} $$
converges absolutely for all $x,y \in \R^{N},$ then according to (\ref{vv}), the series
$$\ds\sum_{\nu\in H}h_{\nu}(x)h_{\nu}(y)(re^{i\alpha})^{|\nu|}$$
converges in $L^{2}(\R^{N}, \ \omega_{k}(x)dx)$ to a function denoted by $K_{\alpha}(r,x,.).$ \\
By the use of Cauchy-Schwartz inequalities, we obtain
\begin{eqnarray*}
D_{k,r}^{\alpha}f(x)&=&
\ds\sum_{\nu \in \Z_{+}^{N}}(re^{i\alpha})^{|\nu|}h_{\nu}(x)\ds\int_{\R^{N}}f(y)h_{\nu}(y)\omega_{k}(y) \ dy  \\&=&\ds\int_{\R^{N}}f(y)\ds\sum_{\nu\in \Z_{+}^{N}}h_{\nu}(x)h_{\nu}(y)(re^{i\alpha})^{|\nu|}\omega_{k}(y) \ dy\\&=& \ds\int_{\R^{N}}K_{\alpha}(r,x,y)f(y)\omega_{k}(y)\ dy.
\end{eqnarray*}
Now, we summarize some properties of the kernel $K_{\alpha}(r,x,y).$
\begin{propo} Let $x, \ y \in \R^{N}, $ $\alpha, \ r \in \R $ such that $\ 0< |\alpha|< \pi$ and $0<r<1,$ then we have\\
$1)$
\begin{eqnarray}
K_{\alpha}(r,x,y)&=&c_{k}\frac{e^{-\frac{1+r^{2}e^{2i\alpha}}{2(1-r^{2}e^{2i\alpha})}(|x|^{2}+|y|^{2})}}{(1-r^{2}e^{2i\alpha})^{\gamma+N/2}}K\left(\frac{2re^{i\alpha}x}{1-r^{2}e^{2i\alpha}},y\right),\label{f}
\end{eqnarray}
$2)$ \begin{eqnarray}
\ds\lim_{r\rightarrow1^{-}}K_{\alpha}(r,x,y)=A_{\alpha}K_{\alpha}(x,y),
\end{eqnarray}
where
\begin{eqnarray}
K_{\alpha}(x,y)&=& e^{-\frac{i}{2}\cot (\alpha)(|x|^{2}+|y|^{2})} K\left(\frac{i
x}{\sin \alpha},y\right), \label{i}
\end{eqnarray}
\begin{eqnarray}
A_{\alpha}=\frac{c_{k}e^{i(\gamma+N/2)(\hat{\alpha}\pi/2-\alpha)}}{(2|\sin
\alpha|)^{\gamma+N/2}} \ \mbox{and} \ \hat{\alpha}=\mbox{sgn} (\sin
\alpha).\label{j}
\end{eqnarray}
$3)$ \begin{eqnarray}
\left|e^{-\frac{(1+r^{2}e^{2i\alpha})|y|^{2}}{2(1-r^{2}e^{2i\alpha})}}K\left(\frac{2re^{i\alpha}x}{1-r^{2}e^{2i\alpha}},y\right)\right|\leq e^{\frac{2r^{2}(1-r^{2})\cos^{2}(\alpha)|x|^{2}}{(r^{4}-2r^{2}\cos 2 \alpha +1)(r^{2}+1)}}.
\end{eqnarray}
\label{pr01}
\end{propo}
\proof
$1)$ According to (\ref{x}), we have
\begin{eqnarray*}
K_{\alpha}(r,x,y)&=&
c_{k}e^{-(|x|^{2}+|y|^{2})/2}\ds\sum_{\nu\in \Z_{+}^{N}}(re^{i\alpha})^{|\nu|}\frac{H_{\nu}(x)H_{\nu}(y)}{2^{|\nu|}}.
\end{eqnarray*}
Using Mehler's formula for the generalized Hermite polynomials (see Theorem 3.12 in \cite{rosler}) which says that for all $x, \ y \in \R^{N}$ and $z\in \C$ with $|z|<1,$
$$\ds\sum_{\nu\in \Z_{+}^{N}}\frac{H_{\nu}(x)H_{\nu}(y)}{2^{|\nu|}}z^{|\nu|}=\frac{e^{-\frac{z^{2}(|x|^{2}+|y|^{2})}{1-z^{2}}}}{(1-z^{2})^{\gamma+(N/2)}}
K\left(\frac{2zx}{1-z^{2}},y\right)$$
and setting $z=re^{i\alpha}$ with $|z|=r<1,$ we obtain the desired result.\\
$2)$ Clearly
\begin{eqnarray*}
\ds\lim_{r\rightarrow 1^{-}}\frac{1+r^{2}e^{2i\alpha}}{1-r^{2}e^{2i\alpha}}&=&i\cot \alpha,\\
\ds\lim_{r\rightarrow 1^{-}}\frac{re^{i\alpha}}{1-r^{2}e^{2i\alpha}}&=&\frac{i}{2\sin \alpha}\\
\ds\lim_{r\rightarrow 1^{-}}(1-r^{2}e^{2i\alpha})^{-(\gamma+\frac{N}{2})}&=&(1-e^{2i\alpha})^{-(\gamma+\frac{N}{2})}\\
&=& \frac{e^{i(\gamma+N/2)(\hat{\alpha}\pi/2-\alpha)}}{(2|\sin
\alpha|)^{\gamma+N/2}}, \ \mbox{where} \ \hat{\alpha}=\mbox{sgn} (\sin \alpha).
\end{eqnarray*}
Then, for $0<|\alpha|<\pi$,
\begin{eqnarray*}
\ds\lim_{r\rightarrow 1^{-}}K_{\alpha}(r,x,y)=A_{\alpha}K_{\alpha}(x,y),
\end{eqnarray*}
where $K_{\alpha}(x,y)$ and $A_{\alpha}$ are defined respectively in (\ref{i}) and (\ref{j}).\\
$3)$ It is straightforward to show that
\begin{eqnarray}
a_{r}=\Re \left(\frac{1+r^{2}e^{2i\alpha}}{1-r^{2}e^{2i\alpha}}\right)&=&\frac{(1-r^{4})}{(1+r^{4})-2r^{2}\cos 2 \alpha}>0,\nonumber\\
b_{r}=\Re\left(\frac{2re^{i\alpha}}{1-r^{2}e^{2i\alpha}}\right)&=&\frac{2(r-r^{3})\cos \alpha }{1+r^{4}-2r^{2}\cos 2 \alpha}. \label{g}
\end{eqnarray}
From (\ref{uu}) and (\ref{g}), we deduce the following majorization:
\begin{eqnarray*}
\left|K\left(\frac{2re^{i\alpha}x}{1-r^{2}e^{2i\alpha}},y\right)\right|\leq e^{|b_{r}| \  |x|  |y|}.
\end{eqnarray*}
Hence,
\begin{eqnarray}
\left|e^{-\frac{(1+r^{2}e^{2i\alpha})|y|^{2}}{2(1-r^{2}e^{2i\alpha})}}K\left(\frac{2re^{i\alpha}x}{1-r^{2}e^{2i\alpha}},y\right)\right|\leq e^{-a_{r}|y|^{2}+|b_{r}| \ |x||y|}.\label{h}
\end{eqnarray}
As $a_{r}>0,$ we deduce that
\begin{eqnarray}
\ds\sup_{s\geq 0} (-a_{r}s^{2}+|b_{r}| \ |x|s)=-\frac{b_{r}^{2}|x|^{2}}{4a_{r}}.\label{k}
\end{eqnarray}
Combining (\ref{h}) and (\ref{k}), we see that
\begin{eqnarray*}
\left|e^{-\frac{(1+r^{2}e^{2i\alpha})|y|^{2}}{2(1-r^{2}e^{2i\alpha})}}K\left(\frac{2re^{i\alpha}x}{1-r^{2}e^{2i\alpha}},y\right)\right|\leq e^{\frac{2r^{2}(1-r^{2})\cos^{2}(\alpha)|x|^{2}}{(r^{4}-2r^{2}\cos 2 \alpha +1)(r^{2}+1)}}.
\end{eqnarray*}
\begin{propo} Let $\alpha \in \R\left\backslash \pi \Z\right.$ and $f \in L^{1}\left(\R^{N},\ \omega_{k}(x)\ dx\right)\cap L^{2}\left(\R^{N},\ \omega_{k}(x)\ dx\right).$ Then the fractional Dunkl transform $D_{k}^{\alpha}$ have the following integral representation
\begin{eqnarray}
D_{k}^{\alpha}f(x)=A_{\alpha}\ds\int_{\R^{N}}f(y)K_{\alpha}(x,y)\omega_{k}(y)dy,a.e,\label{y}
\end{eqnarray}
where
\begin{eqnarray*}
K_{\alpha}(x,y)&=& e^{-\frac{i}{2}\cot (\alpha)(|x|^{2}+|y|^{2})} K\left(\frac{i
x}{\sin \alpha},y\right)
\end{eqnarray*}
and
\begin{eqnarray*}
A_{\alpha}=\frac{c_{k}e^{i(\gamma+N/2)(\hat{\alpha}\pi/2-\alpha)}}{(2|\sin
\alpha|)^{\gamma+N/2}}.
\end{eqnarray*}\label{pro5}
\end{propo}
\proof
$D_{k}^{\alpha}$ is periodic in $\alpha $ with period $2\pi,$ we can assume that $0< |\alpha| < \pi.$ Let $f \in L^{1}\left(\R^{N},\ \omega_{k}(x)\ dx\right)\cap L^{2}\left(\R^{N},\ \omega_{k}(x)\ dx\right).$ From Corollary \ref{coro1},
\begin{eqnarray*}
D_{k}^{\alpha}f(x)=\ds\lim_{j\rightarrow \infty}\ds\int_{\R^{N}}K_{\alpha}(r_{j},x,y)f(y)\omega_{k}(y)\ dy , a.e.
\end{eqnarray*}
From Proposition \ref{pr01}, $2)$ we see that
\begin{eqnarray*}
\ds\lim_{j\rightarrow \infty}K_{\alpha}(r_{j},x,y)f(y)=A_{\alpha}K_{\alpha}(x,y)f(y).
\end{eqnarray*}
Using again Proposition \ref{pr01}, $3)$, we obtain
\begin{eqnarray*}
\left|e^{-\frac{(1+r_{j}^{2}e^{2i\alpha})|y|^{2}}{2(1-r_{j}^{2}e^{2i\alpha})}}K\left(\frac{2r_{j}e^{i\alpha}x}{1-r_{j}^{2}e^{2i\alpha}},y\right)f(y)\right|\leq M_{x} \ |f(y)|,
\end{eqnarray*}
where $M_{x}=\ds\sup_{0\leq r<1}e^{\frac{2r^{2}(1-r^{2})\cos^{2}(\alpha)|x|^{2}}{(r^{4}-2r^{2}\cos 2 \alpha +1)(r^{2}+1)}}.$
\\
Hence, the dominated convergence theorem gives
\begin{eqnarray*}
D_{k}^{\alpha}f(x)=A_{\alpha}\ds\int_{\R^{N}}f(y)K_{\alpha}(x,y)\omega_{k}(y)dy, a.e.
\end{eqnarray*}
\begin{definition} We define the fractional Dunkl transform $D_{k}^{\alpha}$ for $f\in L^{1}\left(\R^{N},\ \omega_{k}(x)\ dx\right)\cap L^{2}\left(\R^{N},\ \omega_{k}(x)\ dx\right)$
by \begin{eqnarray*}
D_{k}^{\alpha}f(x)=A_{\alpha}\ds\int_{\R^{N}}f(y)K_{\alpha}(x,y)\omega_{k}(y)dy.
\end{eqnarray*}
\end{definition}
\begin{remark} \mbox{}\\ $\bullet$ For $\alpha=-\frac{\pi}{2},$ the fractional Dunkl transform $D_{k}^{\alpha}$ is reduces to the Dunkl transform $D_{k}$ and when the multiplicity function $k\equiv 0 , \ D_{k}^{\alpha}$ coincides with the fractional Fourier transform $\mathcal{F}^{\alpha}$ \cite{Brackx}
\begin{eqnarray*}
\mathcal{F}^{\alpha}f(x)=\frac{e^{(iN/2)(\hat{\alpha}\pi/2-\alpha)}}{(2\pi|\sin \alpha|)^{N/2}}\ds\int_{\R^{N}}e^{-\frac{i}{2}(|x|^{2}+|y|^{2})\cot \alpha
+\frac{i }{\sin \alpha}\langle x,y\rangle}f(y) \ dy.
\end{eqnarray*}
$\bullet$ In the one-dimensional case $(N=1),$ the corresponding reflection group $W$ is $\Z_{2}$ and the multiplicity function $k$ is equal to $\nu+1/2>0$. The kernel $K_{\alpha}(x,y)$ defined by (\ref{i}) becomes
\begin{eqnarray}
K_{\alpha}(x,y)&=&e^{-\frac{i}{2}\cot \alpha(x^{2}+y^{2})}E_{\nu}\left(\frac{ix}{\sin \alpha},y\right),\label{zz}
\end{eqnarray}
where $E_{\nu}(x,y)$ is the Dunkl kernel of type $A_{2}$ given by  (see \cite{rosler1})
\begin{eqnarray*}
K(ix,y)=j_{\nu}(xy)+\frac{ixy}{2(\nu+1)}j_{\nu+1}(xy),
\end{eqnarray*}
and $j_{\nu}$ denotes the normalized spherical Bessel function
\begin{eqnarray*}
j_{\nu}(x):=2^{\nu}\Gamma(\nu+1)\ds\frac{J_{\nu}(x)}{x^{\nu}}=\Gamma(\nu+1)\ds
\sum _{n=0}^{+\infty}\ds\frac{(-1)^{n}(x/2)^{2n}}{n!\Gamma(n+\nu+1)}.
\end{eqnarray*}
Here $J_{\nu}$ is the classical Bessel function (see, Watson \cite{Watson}).
The related fractional Dunkl transform $D_{k}^{\alpha}$ in rank-one case takes the form
\begin{eqnarray}
D_{\nu}^{\alpha}f(x)=\mathcal{B}_{\nu}\ds\int_{-\infty}^{+\infty}K_{\alpha}(x,y)f(y)|y|^{2\nu+1} \ dy ,\label{w}
\end{eqnarray}
where
\begin{eqnarray}
\mathcal{B}_{ \nu}=\ds\frac{e^{i(\nu+1)(\hat{\alpha}\pi/2-\alpha)}}{\Gamma(\nu+1)(2|\sin(\alpha)|)^{\nu+1}}.\label{ww}
\end{eqnarray}
Note that if $f$ is an even function then, the fractional Dunkl transform (\ref{w}) coincides with the fractional Hankel transform \cite{Namias1}
\begin{eqnarray*}
H_{\nu}^{\alpha}f(x)=2\mathcal{B}_{\nu} \ds\int_{0}^{+\infty}\ e^{-\frac{i}{2}(x^2+y^{2})
\cot \alpha} \ j_{\nu}\left(\frac{xy}{\sin{\alpha}}\right)f(y)y^{2\nu+1} \ dy.
\end{eqnarray*}
$\bullet$ More generally, for $W=\Z_{2}\times \dots  \times \Z_{2}$ and the multiplicity function $k=(\nu_{1},\dots,\nu_{N}),$ the kernel $K_{\alpha}(x,y)$ defined by (\ref{i})
is given explicitly by
\begin{eqnarray*}
K_{\alpha}(x,y)=e^{-\frac{i}{2}\cot \alpha
(|x|^{2}+|y|^{2})}\ds\prod_{j=1}^{N}E_{\nu_{j}}\left(\frac{i x_{j}}{\sin \alpha},y_{j}\right),
\end{eqnarray*}
where $x=(x_{1},\dots,x_{N}),  \ y=(y_{1},\dots,y_{N}) \in \R^{N}$ and $ E_{\nu_{j}}(x_{j},y_{j})$ is the function defined by (\ref{zz}). In this case the fractional Dunkl transform will be
\begin{eqnarray*}
D_{k}^{\alpha}f(x)=A_{\alpha}\ds\int_{\R^{N}}f(y)K_{\alpha}(x,y) \omega_{k}(y)\ dy,
\end{eqnarray*}
where
\begin{eqnarray*}
A_{\alpha}=\ds\frac{e^{i(\gamma+N/2)(\hat{\alpha}\pi/2-\alpha)}}{\Gamma(\nu_{1}+1)\dots \Gamma(\nu_{N}+1) (2|\sin
\alpha|)^{\gamma+N/2}}
\end{eqnarray*}
and
\begin{eqnarray*}
\omega_{k}(y)=\ds\prod_{j=1}^{N}|x_{j}|^{2\nu_{j}}.
\end{eqnarray*}
\end{remark}
\subsection{Bochner type identity for the fractional Dunkl transform.}
In this subsection, we start with a brief summary on the theory
of k-spherical harmonics. An introduction to
this subject can be found in the monograph \cite{DunklXu}. The space of $k$-spherical
harmonics of degree $n\geq 0$ is defined by
 $$\mathcal{H}_{n}^{k}=\mbox{Ker}\Delta_{k}\cap \mathcal{P}_{n}.$$
Let $S^{N-1}=\left\{x\in \R^{N}; \ |x|=1\right\}$ be the unit sphere in $\R^{N}$ with normalized Lebesgue surface mesure $d\sigma$ and $L^{2}(S^{N-1}, \ \omega_{k}(x) \ d\sigma (x))$ be the Hilbert space with the following inner product given by
$$\langle f,g\rangle_{k}=\ds\int_{S^{N-1}}f(\omega)\overline{g(\omega)}\omega_{k}(\omega) \ d\sigma (\omega).$$
As in the theory of ordinary spherical harmonics, the space $L^{2}(S^{N-1}, \ \omega_{k}(x) \ d\sigma (x))$ decomposes as an orthogonal Hilbert space sum
$$L^{2}(S^{N-1}, \ \omega_{k}(x) \ d\sigma (x))=\ds\bigoplus_{n=0}^{\infty}\mathcal{H}_{n}^{k}.$$
In \cite{yuanxu}, Y. Xu gives an analogue of the Funk-Hecke formula for $k$-spherical harmonics. The well-known special case of the Dunkl-type Funk-Hecke formula is the following (see \cite{rosler2}):
\begin{propo}
Let $N\geq 2$ and put $\lambda=\gamma+(N/2)-1.$ Then for all $Y\in \mathcal{H}_{n}^{k}$ and $x\in \R^{N},$
\begin{eqnarray}
\frac{1}{d_{k}}\ds\int_{S^{N-1}}K(ix,y)Y(y)\omega_{k}(y)\ d\sigma(y)=\frac{\Gamma(\lambda+1)}{2^{n}\Gamma(n+\lambda+1)}j_{n+\lambda}(|x|)Y(ix), \label{hh}
\end{eqnarray}
where
\begin{eqnarray*}
d_{k}=\ds\int_{S^{N-1}}\omega_{k}(y)\ d\sigma(y).
\end{eqnarray*}
In particular
\begin{eqnarray}
\frac{1}{d_{k}}\ds\int_{S^{N-1}}K(ix,y)\omega_{k}(y)\ d\sigma(y)=j_{\lambda}(|x|).\label{gg}
\end{eqnarray}\label{pro6}
\end{propo}
An application of the Dunkl-type Funk-Hecke formula is the following:
\begin{thm} (Bochner type identity) If $f \in L^{1}\left(\R^{N},\ \omega_{k}(x)\ dx\right)\cap L^{2}\left(\R^{N},\ \omega_{k}(x)\ dx\right)$ is of the form $f(x)=p(x)\psi(|x|)$ for some $p\in \mathcal{H}_{n}^{k}$ and a one-variable $\psi$ on $\R_{+},$ then
\begin{eqnarray}
D_{k}^{\alpha}f(x)=e^{in\alpha}p(x)H_{n+\gamma+(N/2)-1}^{\alpha}\psi(|x|).
\end{eqnarray}
In particular, if $f$ is radial, then $$D_{k}^{\alpha}f(x)=H_{n+\gamma+(N/2)-1}^{\alpha}\psi(|x|).$$ \label{th1}
\end{thm}
\proof Since $D_{k}^{\alpha}$ is periodic in $\alpha $ with period $2\pi,$ we can assume that $-\pi< \alpha\leq \pi.$ We see immediately that
\begin{eqnarray*}
D_{k}^{0}f(x)&=&f(x),\\
D_{k}^{\pi}f(x)&=&f(-x)\\&=&p(-x)\psi(-x)\\&=&(-1)^{n}p(x)\psi(x).
\end{eqnarray*}
Now, let $0<|\alpha|< \pi.$ By spherical polar coordinates, we have
\begin{eqnarray}
D_{k}^{\alpha}f(x)&=&A_{\alpha}\ds\int_{\R^{N}}f(y)K_{\alpha}(x,y)\omega_{k}(y)dy\nonumber\\
&=&A_{\alpha}\ds\int_{0}^{+\infty}r^{N-1}F(r,x) \ dr,\label{mm}
\end{eqnarray}
where
\begin{eqnarray*}
F(r,x)&=&\frac{2\pi^{N/2}}{\Gamma(N/2)}\int_{S^{N-1}}K_{\alpha}(x,ry)p(ry)\psi(r|y|)\omega_{k}(ry)\ d\sigma(y).
\end{eqnarray*}
From (\ref{i}) and the homogeneity of $\omega_{k}$ and $p,$ we obtain
\begin{eqnarray*}
F(r,x)&=& \frac{2\pi^{N/2}}{\Gamma(N/2)}e^{-\frac{i}{2}(|x|^{2}+r^{2})\cot (\alpha)} \psi(r)r^{2\gamma+n}\int_{S^{N-1}}p(y)K\left(\frac{irx}{\sin(\alpha)}, y\right)\omega_{k}(y)\ d\sigma(y).
\end{eqnarray*}
Using (\ref{hh}), we get
\begin{eqnarray*}
F(r,x)&=& \frac{2\pi^{N/2}d_{k}}{\Gamma(N/2)}\frac{\Gamma(\lambda+1)}{2^{n}\Gamma(\lambda+n+1)}\\&\times&
e^{-\frac{i}{2}(|x|^{2}+r^{2})\cot (\alpha)} \psi(r)r^{2\gamma+n}p\left(\frac{irx}{\sin(\alpha)}\right)j_{\lambda+n}\left(\frac{r|x|}{\sin(\alpha)}\right),
\end{eqnarray*}
where
\begin{eqnarray*}
\lambda=\gamma+(N/2)-1.
\end{eqnarray*}
Using again the homogeneity of $p,$ we get
\begin{eqnarray*}
F(r,x)
&=&\frac{2\pi^{N/2}d_{k}}{\Gamma(N/2)}\frac{\Gamma(\lambda+1)}{2^{n}\Gamma(\lambda+n+1)}\left(\frac{i}{2\sin \alpha}\right)^{n}\\&\times&
e^{-\frac{i}{2}(|x|^{2}+r^{2})\cot (\alpha)} \psi(r)r^{2\gamma+2n}p(x)j_{\lambda+n}\left(\frac{r|x|}{\sin(\alpha)}\right).
\end{eqnarray*}
Now we can express a relationship between $d_{k}$ and $c_{k}.$ In fact
\begin{eqnarray}
c_{k}^{-1}&=&\ds\int_{\R^{N}}e^{-|y|^{2}}\omega_{k}(y) \ dy\nonumber\\
&=&\frac{2\pi^{N/2}}{\Gamma(N/2)}\ds\int_{0}^{+\infty}r^{N-1}e^{-r^{2}} \ds\int_{S^{N-1}}\omega_{k}(ry)\ d\sigma(y) \ dr \nonumber\\ &=&\frac{2\pi^{N/2}}{\Gamma(N/2)}\ds\int_{0}^{+\infty}r^{2\gamma+N-1}e^{-r^{2}} \ds\int_{S^{N-1}}\omega_{k}(y)\ d\sigma(y) \ dr\nonumber\\&=&\frac{\pi^{N/2}\Gamma(\lambda+1)d_{k}}{\Gamma(N/2)}.\label{ii}
\end{eqnarray}
Recall that
\begin{eqnarray*}
A_{\alpha}=c_{k}\left(\frac{ie^{-i\alpha}}{2 \sin \alpha}\right)^{\gamma+(N/2)},
\end{eqnarray*}
then by the use of (\ref{ii}), we obtain
\begin{eqnarray*}
A_{\alpha}\frac{2\pi^{N/2}d_{k}}{\Gamma(N/2)}\frac{\Gamma(\lambda+1)}{2^{n}\Gamma(\lambda+n+1)} \left(\frac{i}{2\sin \alpha}\right)^{n}&=&\frac{2\left(\frac{ie^{-i\alpha}}{2\sin \alpha}\right)^{\lambda+n+1}}{\Gamma(\lambda+n+1)}e^{in\alpha}\\&=& 2\mathcal{B}_{\nu} e^{in\alpha}.
\end{eqnarray*}
Hence
\begin{eqnarray}
F(r,x)=2\mathcal{B}_{\nu} e^{in\alpha}e^{-\frac{i}{2}(|x|^{2}+r^{2})\cot (\alpha)} \psi(r)r^{2\gamma+2n}p(x)j_{\lambda+n}\left(\frac{r|x|}{\sin(\alpha)}\right).\label{ll}
\end{eqnarray}
Substituting (\ref{ll}) in (\ref{mm}) to get
\begin{eqnarray*}
D_{k}^{\alpha}f(x)&=&2\mathcal{B}_{\nu}e^{in\alpha}p(x)\\&\times&\ds\int_{0}^{+\infty}e^{-\frac{i}{2}(|x|^{2}+r^{2})\cot (\alpha)}\psi(r)r^{2(\lambda+n)+1} j_{\lambda+n}\left(\frac{r|x|}{\sin(\alpha)}\right)\ dr\nonumber\\
& =&e^{in\alpha}p(x)H_{n+\lambda}^{\alpha}\psi(|x|)\\
& =&e^{in\alpha}p(x)H_{n+\gamma+(N/2)-1}^{\alpha}\psi(|x|).
\end{eqnarray*}

Now, we give the material needed for an application of Bochner type identity. Let $\{p_{n,j}\}_{j\in J_{n}}$ be an orthonormal
basis of $ \mathcal{H}_{n}^{k}.$ Let $m, \ n$ be non-negative integers and $j\in J_{n}.$ Define
\begin{eqnarray*}
c_{m,n}=\left(\frac{m! \ \Gamma(N/2)}{\pi^{N/2}\Gamma((N/2)+\gamma+n+m)}\right)^{1/2}
\end{eqnarray*}
and
\begin{eqnarray}
\psi_{m,n,j}(x)=c_{m,n} \ p_{n,j}(x) \ L_{m}^{(n+\gamma+N/2-1)}(|x|^{2}) \ e^{-|x|^{2}/2},
\end{eqnarray}
where $L_{m}^{(a)}$ denote the Laguerre polynomial defined by
\begin{eqnarray*}
L_{m}^{(a)}(x)=\frac{x^{-a}e^{x}}{n!}\frac{d^{n}}{dx^{n}}\left(e^{-x}x^{n+a}\right).
\end{eqnarray*}
It follows from Proposition 2.4 and Theorem 2.5 of Dunkl \cite{Dunkl3} that
$$\left\{\psi_{m,n,j}:m, \ n=0, \ 1, \ 2, \dots, \ j \in J_{n}\right\}$$
forms an orthonormal basis of $ L^{2}(\R^{N}, \ \omega_{k}(x)dx).$
\begin{thm} The family $\left\{\psi_{m,n,j}:m, \ n=0, \ 1, \ 2, \dots, \ j \in J_{n}\right\}$ are a basis of eigenfunctions of the fractional Dunkl transform $D_{k}^{\alpha}$ on $L^{2}\left(\R^{N},\ \omega_{k}(x)\ dx\right),$ satisfying
\begin{eqnarray}
D_{k}^{\alpha}\psi_{m,n,j}=e^{i\alpha(n+2m)}\psi_{m,n,j}.\label{jj}
\end{eqnarray}
\end{thm}
\proof We need only to prove (\ref{jj}). Applying Theorem \ref{th1} with $p$ replaced by $p_{n,j}$
and with $\psi(r)=L_{m}^{(n+\gamma+N/2-1)}(r^{2}) \ e^{-r^{2}/2}$, we obtain
\begin{eqnarray*}
 D_{k}^{\alpha}\psi_{m,n,j}(x)=c_{m,n}e^{in\alpha}p_{n,j}(x)H_{\nu}^{\alpha}\psi(|x|),
\end{eqnarray*}
where
\begin{eqnarray*}
\nu=n+\gamma+(N/2)-1,
\end{eqnarray*}
and
\begin{eqnarray*}
 \mathcal{H}_{\nu}^{\alpha}\psi(|x|)&=&2\mathcal{B}_{\nu}\ds\int_{0}^{+\infty}e^{-\frac{i}{2}\cot(\alpha)(|x|^{2}+r^{2})}
 j_{\nu}\left(\frac{r|x|}{\sin \alpha}\right)L_{m}^{(\nu)}(r^{2})e^{-\frac{r^{2}}{2}}r^{2\nu+1} \ dr.
\end{eqnarray*}
Observe that
\begin{eqnarray*}
 \mathcal{H}_{\nu}^{\alpha}\psi(|x|)&=&2\mathcal{B}_{\nu}e^{-\frac{i}{2}\cot(\alpha)|x|^{2}}I_{\nu},
\end{eqnarray*}
where
\begin{eqnarray*}
I_{\nu}&=&\ds\int_{0}^{+\infty}r^{2\nu+1}L_{m}^{(\nu)}(r^{2})e^{-(\frac{1}{2}+\frac{i}{2}\cot(\alpha))r^{2}} j_{\nu}\left(\frac{r|x|}{\sin \alpha}\right)\ dr\\
&=& 2^{\nu}\Gamma(\nu+1)\left(\frac{\sin \alpha}{|x|}\right)^{\nu}\ds\int_{0}^{+\infty}r^{\nu+1}L_{m}^{(\nu)}(r^{2})e^{-(\frac{1}{2}+\frac{i}{2}\cot(\alpha))r^{2}} J_{\nu}\left(\frac{r|x|}{\sin \alpha}\right)\ dr.
\end{eqnarray*}
To compute $I_{\nu}$, we need the following formulas (see 7.4.21 (4) in \cite{Gradshteyn})   \begin{eqnarray*}
\ds\int_{0}^{+\infty}y^{\nu+1}e^{-\beta y^{2}}L_{m}^{\nu}(ay^{2})J_{\nu}(zy)\ dy=d_{m}z^{\nu}e^{-z^{2}/(4\beta)}L_{m}^{\nu}\left[\frac{az^{2}}{4\beta(a-\beta)}\right]
\end{eqnarray*}
where $d_{m}=((\beta-a)^{m}/(2^{\nu+1}\beta^{\nu+m+1})), \ a, \ \Re \beta>0,\ \Re \nu>-1. $ \\
Let us take $\beta=\frac{1}{2}+\frac{i}{2}\cot(\alpha)=\frac{ie^{-i\alpha}}{2\sin \alpha}, \ a=1$ and $z=\frac{|x|}{\sin \alpha},$ then
\begin{eqnarray*}
d_{m}&=&\frac{e^{2i\alpha m}}{2^{\nu+1}\mathcal{A}_{\alpha}\Gamma(\nu+1)},\\
\frac{az^{2}}{4\beta(a-\beta)}&=&|x|^{2},\\
-\frac{z^{2}}{4\beta}&=&-\frac{|x|^{2}}{2}+\frac{i}{2}\cot(\alpha)|x|^{2}.
\end{eqnarray*}
Hence
\begin{eqnarray*}
\ds\int_{0}^{+\infty}r^{\nu+1}L_{m}^{(\nu)}(r^{2})e^{-(\frac{1}{2}+\frac{i}{2}\cot(\alpha))r^{2}} J_{\nu}\left(\frac{r|x|}{\sin \alpha}\right)\ dr=\frac{e^{2i\alpha m}e^{-\frac{|x|^{2}}{2}+\frac{i}{2}\cot(\alpha)|x|^{2}}}{2^{\nu+1}\mathcal{A}_{\alpha}\Gamma(\nu+1)}
\left(\frac{|x|}{\sin \alpha}\right)^{\nu}L_{m}^{(\nu)}(|x|^{2}),
\end{eqnarray*}
and therefore
\begin{eqnarray*}
 H_{\nu}^{\alpha}\psi(|x|)&=&e^{2i\alpha m}L_{m}^{(\nu)}(|x|^{2}) \ e^{-|x|^{2}/2},
\end{eqnarray*}
which finishes the proof.
\subsection{Master Formula for the fractional Dunkl transform.}
In this subsection, we are interesting with a master formula for the fractional Dunkl transform. For this we need the following lemma
\begin{lemma} Let $p\in \mathcal{P}_{n}$ and $x=(x_{1},\dots,x_{N})\in \C^{N}.$ Then for $\omega\in \C$ and $\Re(\omega)>0,$
\begin{eqnarray}
c_{k}\ds\int_{\R^{N}}p(y)K(x,2y)e^{-\omega|y|^{2}}\omega_{k}(y) \ dy =\ds\frac{e^{\frac{l(x)}{\omega}}}{\omega^{\gamma+n+(N/2)}}e^{\frac{\omega}{4}\Delta_{k}}p(x),\label{r}
\end{eqnarray}
where $l(x)=\ds\sum_{j=1}^{N}x_{j}.$
\label{lemma2}
\end{lemma}
\proof
First compute the above integral when $\omega>0.$
\begin{eqnarray*}
 c_{k}\ds\int_{\R^{N}}p(y)K(x,2y)e^{-\omega|y|^{2}}\omega_{k}(y) \ dy =c_{k}\ds\int_{\R^{N}}p(y)K(x,2y)e^{-|\sqrt{\omega}y|^{2}}\omega_{k}(y) \ dy.
\end{eqnarray*}
By the change of variables $u=\sqrt{\omega}y$ and the homogeneity of $\omega_{k}$ and $p,$ we obtain
\begin{eqnarray}
& \ & c_{k}\ds\int_{\R^{N}}p(y)K(x,2y)e^{-\omega|y|^{2}}\omega_{k}(y) \ dy  \nonumber\\
&=&\frac{c_{k}}{\omega^{\gamma+(n+N)/2}}
\ds\int_{\R^{N}}p(y)K\left(\frac{x}{\sqrt{\omega}},2y\right)e^{-|y|^{2}}\omega_{k}(y)
\ dy.\label{o}
\end{eqnarray}
Using Proposition 2.1 from \cite{Dunkl3}, which says that for all $p\in \mathcal{P}$ and $x\in \C^{N}$
\begin{eqnarray*}
c_{k}\ds\int_{\R^{N}}e^{-\frac{\Delta_{k}}{4}}p(y)K(x,2y)e^{-|y|^{2}}\omega_{k}(y) \ dy=e^{l(x)}p(x),
\end{eqnarray*}
we deduce the following identity:
\begin{eqnarray}
c_{k}\ds\int_{\R^{N}}p(y)K(x,2y)e^{-|y|^{2}}\omega_{k}(y) \ dy=e^{l(x)}e^{\frac{\Delta_{k}}{4}}p(x).\label{p}
\end{eqnarray}
Combining (\ref{o}) and (\ref{p}) to get
\begin{eqnarray*}
c_{k}\ds\int_{\R^{N}}p(y)K(x,2y)e^{-\omega|y|^{2}}\omega_{k}(y) \ dy=\ds\frac{e^{\frac{l(x)}{\omega}}}{\omega^{\gamma+(n+N)/2}}e^{\frac{\Delta_{k}}{4}}p\left(\frac{x}{\sqrt{\omega}}\right).
\end{eqnarray*}
Now use Lemma 2.1 from \cite{rosler} to obtain
$$e^{\frac{\Delta_{k}}{4}}p\left(\frac{x}{\sqrt{\omega}}\right)=\frac{1}{\omega^{n/2}}e^{\frac{\omega}{4}\Delta_{k}}p(x).$$
Hence, we find the equality (\ref{r}) for $\omega>0$. By analytic continuation, this holds for $\{\omega\in\C: \ \Re(\omega)>0\}.$

We are now in a position to give the master formula.
\begin{thm} Let $p\in \mathcal{P}_{n}$ and $x\in \R^{N}.$ Then
\begin{eqnarray}
D_{k}^{\alpha}\left[e^{-\frac{|y|^{2}}{2}}e^{-\frac{\Delta_{k}}{4}}p(y)\right](x)=e^{in\alpha}e^{-\frac{|x|^{2}}{2}}e^{-\frac{\Delta_{k}}{4}}p(x). \label{mas1}
\end{eqnarray}\label{th1}
\end{thm}
\proof
It follows easily from  (\ref{y}) that
\begin{eqnarray*}
D_{k}^{\alpha}\left[e^{-\frac{|y|^{2}}{2}}e^{-\frac{\Delta_{k}}{4}}p(y)\right](x)=A_{\alpha}e^{-\frac{i}{2}\cot(\alpha)|x|^{2}}
\ds\int_{\R^{N}}e^{-\frac{\Delta_{k}}{4}}p(y)K\left(\frac{ix}{\sin \alpha},y\right)e^{-\omega |y|^{2}}\omega_{k}(y) \ dy ,
\end{eqnarray*}
where
\begin{eqnarray}
\omega=\frac{1}{2}+\frac{i}{2}\cot(\alpha)=\frac{ie^{-i\alpha}}{2\sin \alpha}. \label{u}
\end{eqnarray}
Since
$$e^{-\frac{\Delta_{k}}{4}}p(y)=\ds\sum_{s=0}^{[\frac{n}{2}]}\frac{(-1)^{s}}{s! 4^{s}}\Delta_{k}^{s}p(y),$$
we conclude that
\begin{eqnarray}
\ds\int_{\R^{N}}e^{-\frac{\Delta_{k}}{4}}p(y)K\left(\frac{ix}{\sin \alpha},y\right)e^{-\omega |y|^{2}}\omega_{k}(y) \ dy=\ds\sum_{s=0}^{[\frac{n}{2}]}\frac{(-1)^{s}}{s! 4^{s}}\ds\int_{\R^{N}}\Delta_{k}^{s}p(y)K\left(\frac{ix}{\sin \alpha},y\right)e^{-\omega |y|^{2}}\omega_{k}(y) \ dy. \nonumber\\ \label{z}
\end{eqnarray}
For $s\in \Z_{+}$ with $2s\leq n,$ the polynomial $\Delta_{k}^{s}p$ is homogeneous of degree $n-2s.$ Hence by the previous Lemma, we obtain
\begin{eqnarray}
c_{k}\ds\int_{\R^{N}}\Delta_{k}^{s}p(y)K\left(\frac{ix}{\sin \alpha},y\right)e^{-\omega |y|^{2}}\omega_{k}(y) \ dy=\ds\frac{e^{\frac{l(X_{\alpha})}{\omega}}}{\omega^{\gamma+n+(N/2)}}e^{\frac{\omega}{4}\Delta_{k}}\left[\omega^{2s} \Delta_{k}^{s} p\right](X_{\alpha}), \label{s}
\end{eqnarray}
where
\begin{eqnarray}
X_{\alpha}=\frac{i x}{2\sin \alpha}. \label{t}
\end{eqnarray}
Substituting (\ref{s}) in (\ref{z}) to get
\begin{eqnarray*}
c_{k}\ds\int_{\R^{N}}e^{-\frac{\Delta_{k}}{4}}p(y)K\left(\frac{ix}{\sin \alpha},y\right)e^{-\omega |y|^{2}}\omega_{k}(y) \ dy&=&\ds\frac{e^{\frac{l(X_{\alpha})}{\omega}}}{\omega^{\gamma+n+(N/2)}}e^{\frac{\omega}{4}\Delta_{k}}\ds\sum_{s=0}^{[\frac{n}{2}]}\frac{(-1)^{s}\omega^{2s}}{s! 4^{s}}\Delta_{k}^{s}p(X_{\alpha})\\&=&\ds\frac{e^{\frac{l(X_{\alpha})}{\omega}}}{\omega^{\gamma+n+(N/2)}}e^{\frac{\omega}{4}\Delta_{k}}e^{-\frac{\omega^{2}}{4}\Delta_{k}}p(X_{\alpha})\\
&=& \ds\frac{e^{\frac{l(X_{\alpha})}{\omega}}}{\omega^{\gamma+n+(N/2)}}e^{\frac{\omega-\omega^{2}}{4}\Delta_{k}}p(X_{\alpha}).
\end{eqnarray*}
Replacing $\omega$ and $X_{\alpha}$ by their values given in (\ref{u}) and (\ref{t}) and use Lemma 2.1 in \cite{rosler}, we obtain
\begin{eqnarray*}
e^{\frac{\omega-\omega^{2}}{4}\Delta_{k}}p(X_{\alpha})&=&\frac{i^{n}}{2^{n}\sin^{n}\alpha} e^{-\sin^{2}(\alpha)(\omega-\omega^{2})\Delta_{k}}p(x)\\
&=&\frac{i^{n}}{2^{n}\sin^{n}\alpha} e^{-\frac{\Delta_{k}}{4}}p(x).
\end{eqnarray*}
Also,
\begin{eqnarray*}
\omega^{n+\gamma+(N/2)}&=&\left(\frac{ie^{-i\alpha}}{2\sin \alpha}\right)^{n+\gamma+(N/2)}=\frac{i^{n}e^{-in\alpha}}{2^{n}\sin^{n}\alpha} \ \frac{e^{i(\gamma+N/2)(\hat{\alpha}\pi/2-\alpha)}}{(2|\sin \alpha|)^{\gamma+(N/2)}}\\
e^{\frac{l(X_{\alpha})}{\omega}}&=& e^{\frac{ie^{i\alpha}}{2\sin \alpha}|x|^{2}}.
\end{eqnarray*}
Then
\begin{eqnarray}
\ds\int_{\R^{N}}e^{-\frac{\Delta_{k}}{4}}p(y)K\left(\frac{ix}{\sin \alpha},y\right)e^{-\omega |y|^{2}}\omega_{k}(y) \ dy=A_{\alpha}^{-1}e^{in\alpha}e^{\frac{ie^{i\alpha}}{2\sin \alpha}|x|^{2}}e^{-\frac{\Delta_{k}}{4}}p(x).\label{aa}
\end{eqnarray}
 Finally, if we multiply equation (\ref{aa}) by $A_{\alpha}e^{-\frac{i}{2}\cot(\alpha)|x|^{2}},$ we obtain the desired result.

A consequence of the Master formula (\ref{mas1}) is
\begin{coro} (Hecke type identity) If in addition to the assumption in Theorem \ref{th1}, the polynomial $p\in \mathcal{H}_{n}^{k},$ then (\ref{mas1}) becomes
\begin{eqnarray}
D_{k}^{\alpha}\left[e^{-\frac{|.|^{2}}{2}}p\right](x)=e^{in\alpha}e^{-\frac{|x|^{2}}{2}}p(x). \label{mas2}
\end{eqnarray}
\end{coro}
Now, we are interesting to complete the spectral study of $T$ started in proposition \ref{pro7} by means of the Master formula. In fact we have the following
\begin{coro} \mbox{}\\ $ L^{2}\left(\R^{N},\ \omega_{k}(x)\ dx\right)$ decomposes as an orthogonal Hilbert space sum according to
$$L^{2}\left(\R^{N},\ \omega_{k}(x)\ dx\right)=\ds\bigoplus_{n\in \Z_{+}}V_{n},$$
where
$$V_{n}=\left\{e^{-\frac{|x|^{2}}{2}}e^{-\frac{\Delta_{k}}{4}}p(x); \quad p \in \mathcal{P}_{n}\right\}$$
is the eigenspace of $T$ corresponding to the eigenvalue $in.$ In particular, $T$
is essentially self-adjoint. The spectrum of its closure is purely discrete
and given by
$$\sigma(\overline{T})=i\Z_{+.}$$
\end{coro}
\proof\mbox{}\\ Let $f$ be an element of the subspace $V_{n}$ defined by
$$f(x)=e^{-\frac{|x|^{2}}{2}}e^{-\frac{\Delta_{k}}{4}}p(x),$$
where $p\in \mathcal{P}_{n}.$ From (\ref{mas1}), the limits
\begin{eqnarray*}
\ds\lim_{\alpha\rightarrow 0}\frac{D_{k}^{\alpha}f-f}{\alpha}&=&\ds\lim_{\alpha\rightarrow 0}\frac{e^{in\alpha }-1}{\alpha}f
\end{eqnarray*}
exists in $L^{2}(\R^{N}, \ \omega_{k}(x)dx)$ and equals $inf.$ Then
\begin{eqnarray}
f\in D(T) \quad \mbox{and} \quad T(f)=inf. \label{y_{1}}
\end{eqnarray}
Hence, $V_{n}$ is the eigenspace of $T$ corresponding to the eigenvalue $in.$
\section{Realization of the operator T.}
The aim of the following is to find a subspace $\mathcal{W}\subset D(T)$ of $L^{2}(\R^{N}, \ \omega_{k}(x)dx)$ in which we define $T$ explicitly.
\begin{lemma} For $z\in\C^{N}$ set $l(z)=\ds\sum_{i=1}^{N}z_{i}^{2}$. Then for all $z,\omega \in \C^{N},$
\begin{equation}
c_{k}\ds\int_{\R^{N}}K(2z,x)K(2\omega,x)e^{-A|x|^{2}}\omega_{k}(x) \
dx=\ds\frac{e ^{\frac{l(z)+l(\omega)}{A}}}{A^{\gamma+N/2}}
 K(2z/A,\omega),
\end{equation}
where $A$ is a complex number such that $\Re(A)>0$.\label{lemma1}
\end{lemma}
\proof The result is obtained by means of a similar technic used in the proof of Lemma \ref{lemma2} and the following formula (see \cite{Dunkl3})
$$c_{k}\ds\int_{\R^{N}}K(2z,x)K(2\omega,x)e^{-|x|^{2}}\omega_{k}(x) \ dx=e^{l(z)+l(\omega)}K(2z,\omega).$$
\begin{thm} Let $f\in L^{1}(\R^{N}, \ \omega_{k}(x)dx)\cap L^{2}(\R^{N}, \ \omega_{k}(x)dx)$ such that $D_{k}f \in L^{1}(\R^{N}, \ \omega_{k}(x)dx)$ and $\alpha \not\in \left\{ \frac{\pi}{2}+k\pi, \ k \in \Z\right\}.$ Then
\begin{eqnarray}
D_{k}^{\alpha}f(x)=c_{k}\left(\frac{e^{-i\alpha}}{2\cos \alpha}\right)^{\gamma+\frac{N}{2}}\ds\int_{\R^{N}}e^{\frac{i}{2}\tan (\alpha)(|x|^{2}+|y|^{2})}K(\frac{ix}{\cos \alpha},y)D_{k}f(y)\omega_{k}(y)dy.\label{tt}
\end{eqnarray}
\label{tt2}
\end{thm}
\proof Let $f\in L^{1}(\R^{N}, \ \omega_{k}(x)dx)\cap L^{2}(\R^{N}, \ \omega_{k}(x)dx)$ such that $D_{k}f \in L^{1}(\R^{N}, \ \omega_{k}(x)dx).$ Let $\epsilon$ be an arbitrary positive number and put
$$F_{\epsilon}(x)=\ds\int_{\R^{N}}f(y)g_{\epsilon}(y)\omega_{k}(y)\ dy ,$$
where $g_{\epsilon}(y)=\ e^{-(\epsilon+\frac{i}{2}\cot
\alpha)|y|^{2}}K\left(\frac{ix}{\sin \alpha},y\right).$ \\
From (\ref{v}), we deduce that $$|f(y)g_{\epsilon}(y)|\leq |f(y)|,$$
so the dominated convergence theorem can be invoked again to give
\begin{eqnarray}
\ds\lim_{\epsilon\rightarrow 0}F_{\epsilon}(x)&=&\frac{e ^{ \frac{i}{2}|x|^{2}\cot \alpha}}{A_{\alpha}}D_{k}^{\alpha}f(x).\label{qq}
\end{eqnarray}
Using Lemma \ref{lemma1}, we can show
\begin{eqnarray}
D_{k}g_{\epsilon}(\xi)
&=&\frac{ e ^{-\frac{|x|^{2}}{4\epsilon\sin^{2}\alpha+i\sin
2 \alpha}}}{(2\epsilon+i\cot
\alpha)^{\gamma+N/2}}\ \ e ^{-\frac{|\xi|^{2}}{4\epsilon+2i\cot
\alpha}}\nonumber\\& \times & K\left(\frac{x}{2\epsilon\sin
\alpha+i\cos \alpha },\xi\right).\label{kk}
\end{eqnarray}
Now applying the Parseval formula for the Dunkl transform (see Lemma 4.25, \cite{deJeu}) and using (\ref{kk}), we obtain
\begin{eqnarray*}
F_{\epsilon}(x)&=&\frac{e^{-\frac{|x|^{2}}{4\epsilon\sin^{2}\alpha+i\sin
2 \alpha}}}{(2\epsilon+i\cot
\alpha)^{\gamma+N/2}} \\\ & \times&
\ds\int_{\R^{N}}\ e ^{-\frac{|\xi|^{2}}{4\epsilon+2i\cot
\alpha}} \ K\left(\frac{x}{2\epsilon\sin \alpha+i\cos \alpha
},\xi\right)D_{k}f(-\xi)\omega_{k}(\xi) \ d\xi.
\end{eqnarray*}
(\ref{uu}) gives again the following majorization:
\begin{eqnarray*}
\left|K\left(\frac{x}{2\epsilon\sin \alpha+i\cos \alpha
},\xi\right)\right|\leq e^{\frac{2\epsilon \sin \alpha}{4\epsilon^{2}\sin^{2}(\alpha)+\cos^{2}(\alpha)}|x||\xi|}.
\end{eqnarray*}
Hence,
\begin{eqnarray}
\left|e ^{-\frac{|\xi|^{2}}{4\epsilon+2i\cot
\alpha}} \ K\left(\frac{x}{2\epsilon\sin \alpha+i\cos \alpha
},\xi\right)\right|\leq e^{-p_{\epsilon}|\xi|^{2}+q_{\xi}|\xi|},\label{pp}
\end{eqnarray}
where
\begin{eqnarray*}
p_{\epsilon}=\frac{\epsilon}{4\epsilon^{2}+\cot^{2}\alpha}
\end{eqnarray*}
and
\begin{eqnarray*}
q_{\epsilon}=\frac{2\epsilon \sin (\alpha)|x|}{4\epsilon^{2}\sin^{2}(\alpha)+\cos^{2}(\alpha)}.
\end{eqnarray*}
As $p_{\epsilon}>0,$ we deduce that
\begin{eqnarray}
\ds\sup_{s\geq 0} (-p_{\epsilon}s^{2}+q_{\epsilon}s)=-\frac{q_{\epsilon}^{2}}{4p_{\epsilon}}.\label{q}
\end{eqnarray}
Applying formula (\ref{pp}) and (\ref{q}), we obtain
\begin{eqnarray*}
\left|e ^{-\frac{|\xi|^{2}}{4\epsilon+2i\cot
\alpha}} \ K\left(\frac{x}{2\epsilon\sin \alpha+i\cos \alpha
},\xi\right)D_{k}f(-\xi)\right|&\leq& e^{-\frac{4\epsilon |x|^{2}}{4\epsilon^{2}\sin^{2}(\alpha)+\cos^{2}(\alpha)}}|D_{k}f(-\xi)|
\\&\leq& B_{x}|D_{k}f(-\xi)|,
\end{eqnarray*}
where $B_{x}=\ds\sup_{\epsilon\in]0,1]}e^{-\frac{4\epsilon |x|^{2}}{4\epsilon^{2}\sin^{2}(\alpha)+\cos^{2}(\alpha)}}.$
The function $\xi \mapsto D_{k}f(-\xi)$ is in $L^{1}(\R^{N}, \ \omega_{k}(x)dx),$ then the dominated convergence theorem implies
\begin{eqnarray}
& \ & \ds\lim_{\epsilon\rightarrow
0}F_{\epsilon}(x)=\frac{e^{\frac{i|x|^{2}}{\sin 2
\alpha}}}{(i\cot
\alpha)^{\gamma+N/2}}\ \nonumber\\ &\times&
\ds\int_{\R^{N}}\ e^{\frac{i |\xi|^{2}\tan
\alpha}{2}}\ K\left(-\frac{ix}{\cos
\alpha},\xi\right)D_{k}f(-\xi)\omega_{k}(\xi) \ d\xi. \label{rr}
\end{eqnarray}
Hence, (\ref{qq}) and (\ref{rr}) gives after simplification
\begin{eqnarray}
D_{k}^{\alpha}f(x)&=&c_{k}\left(\frac{e^{-i\alpha}}{2\cos \alpha}\right)^{\gamma+\frac{N}{2}}\ e^{\frac{i}{2}|x|^{2}\tan
\alpha}\nonumber\\
&\times&\ds\int_{\R^{N}}\ e ^{\frac{i}{2}|\xi|^{2}\tan
\alpha} \ K\left(-\frac{ix}{\cos
\alpha},\xi\right)D_{k}f(-\xi)\omega_{k}(\xi)\ d\xi.\label{ss}
\end{eqnarray}
Finally, if we make the change of variables $u=-y$ in ($\ref{ss}$), then we find (\ref{tt}).
\begin{remark} Using (\ref{tt}) together with the dominated convergence theorem, we get
\begin{eqnarray*}
\ds\lim _{\alpha\rightarrow 0^{+}}D_{k}^{\alpha}f(x)&=&\ds\lim _{\alpha\rightarrow 0^{-}}D_{k}^{\alpha}f(x)=D_{k}^{2}f(-x)=f(x), a.e,\\
\ds\lim _{\alpha\rightarrow \pi^{-}}D_{k}^{\alpha}f(x)&=&\ds\lim _{\alpha\rightarrow -\pi^{+}}D_{k}^{\alpha}f(x)=D_{k}^{2}f(x)=f(-x) \ a.e.
\end{eqnarray*}
\end{remark}
\begin{coro} Under the assumptions of Theorem \ref{tt2}, we have
\begin{eqnarray}
\frac{D_{k}^{\alpha}f(x)-f(x)}{\alpha}&=&r_{1}(\alpha) \ \frac{c_{k}}{2^{\gamma+(N/2)}}
\ds\int_{\R^{N}}e^{\frac{i}{2}\tan (\alpha)(|x|^{2}+|y|^{2})}K\left(\frac{ix}{\cos \alpha},y\right)D_{k}f(y)\omega_{k}(y)dy\nonumber\\&+&\frac{c_{k}}{2^{\gamma+\frac{N}{2}}}
\ds\int_{\R^{N}}r_{2}(\alpha,x,y)D_{k}f(y)\omega_{k}(y)dy, a.e,
\end{eqnarray}
where
\begin{eqnarray*}
r_{1}(\alpha)=\frac{\left(\frac{e^{-i\alpha}}{\cos \alpha}\right)^{\gamma+\frac{N}{2}}-1}{\alpha}
\end{eqnarray*}
and
\begin{eqnarray*}
r_{2}(\alpha,x,y)=\frac{e^{\frac{i}{2}\tan (\alpha)(|x|^{2}+|y|^{2})}K(\frac{ix}{\cos \alpha},y)-K(ix,y)}{\alpha}.
\end{eqnarray*}
\end{coro}
\proof The result is consequence of (\ref{tt}) and (\ref{w3}).
\begin{lemma} Let $\alpha_{0}\in]0,\frac{\pi}{2}[$ and $x,y \in \R^{N}.$ Then
\begin{eqnarray}
\left|r_{2}(\alpha,x,y)\right|\leq \frac{1}{2}(1+\tan ^{2}\alpha_{0})(|x|^{2}+|y|^{2})+\frac{|\sin (\alpha_{0})|}{\cos^{2} (\alpha_{0})} \ \sqrt{N} \ |x||y|,
\end{eqnarray}
where $\alpha\in ]0,\alpha_{0}].$
\end{lemma}
\proof By the mean value theorem, we have
\begin{eqnarray*}
\left|r_{2}(\alpha,x,y)\right|\leq \ds\sup_{\alpha \in [0,\alpha_{0}]}\left|\frac{\partial}{\partial\alpha}r_{3}(\alpha,x,y)\right|,
\end{eqnarray*}
where
\begin{eqnarray*}
r_{3}(\alpha,x,y)=e^{\frac{i}{2}\tan (\alpha)(|x|^{2}+|y|^{2})}K\left(\frac{ix}{\cos \alpha},y\right).
\end{eqnarray*}
From (\ref{xx}), we get
\begin{eqnarray*}
K\left(\frac{ix}{\cos \alpha},y\right)=K\left(x,\frac{iy}{\cos \alpha}\right).
\end{eqnarray*}
Therefore,
\begin{eqnarray*}
r_{3}(\alpha,x,y)=e^{\frac{i}{2}\tan (\alpha)(|x|^{2}+|y|^{2})}K\left(x,\frac{iy}{\cos \alpha}\right).
\end{eqnarray*}
A simple calculations shows that
\begin{eqnarray}
\frac{\partial}{\partial \alpha}r_{3}(\alpha,x,y)&=&\frac{i}{2}(1+\tan ^{2}\alpha)(|x|^{2}+|y|^{2})r_{3}(\alpha,x,y)\nonumber\\&+& \frac{i\sin(\alpha)}{\cos^{2}(\alpha)} \ e^{\frac{i}{2}\tan (\alpha)(|x|^{2}+|y|^{2})} \ds\sum_{j=1}^{N}y_{j}\frac{\partial}{\partial y_{j}}K\left(x,\frac{iy}{\cos \alpha}\right).\label{dd}
\end{eqnarray}
From (\ref{uu}), the inequality
\begin{eqnarray*}
\left|\frac{\partial}{\partial y_{j}}K\left(x,\frac{iy}{\cos \alpha}\right)\right|\leq |x|
\end{eqnarray*}
holds and hence
\begin{eqnarray*}
\left|\frac{\partial}{\partial \alpha}r_{3}(\alpha,x,y)\right|&\leq& \frac{1}{2}(1+\tan ^{2}\alpha)(|x|^{2}+|y|^{2})+\frac{|\sin (\alpha)|}{\cos^{2} (\alpha)} \ |x|\ds\sum_{j=1}^{N}|y_{j}|\\
&\leq& \frac{1}{2}(1+\tan ^{2}\alpha)(|x|^{2}+|y|^{2})+\frac{|\sin (\alpha)|}{\cos^{2} (\alpha)} \ \sqrt{N} \ |x||y| \\
&\leq& \frac{1}{2}(1+\tan ^{2}\alpha_{0})(|x|^{2}+|y|^{2})+\frac{|\sin (\alpha_{0})|}{\cos^{2} (\alpha_{0})} \ \sqrt{N} \ |x||y|.
\end{eqnarray*}
Which finishes the proof.
\begin{thm} Let
\begin{eqnarray*}
\mathcal{W}=\left\{f\in L^{1}(\R^{N}, \ \omega_{k}(x)dx)\cap L^{2}(\R^{N}, \ \omega_{k}(x)dx \ ; \ |y|^{2}f \in L^{2}(\R^{N}, \ \omega_{k}(x)dx)\right. \\ \left. \mbox{and} \  |y|^{2}D_{k}f \in L^{1}(\R^{N}, \ \omega_{k}(x)dx)\cap L^{2}(\R^{N}, \ \omega_{k}(x)dx)\right\}.
\end{eqnarray*}
Then for all $f\in \mathcal{W},$
\begin{eqnarray}
Tf(x)=-i(\gamma+(N/2))f(x)+\frac{i}{2}|x|^{2}f(x)+\frac{i}{2}D_{k}\left[|y|^{2}D_{k}(y)\right](-x) \ a.e.\label{ff}
\end{eqnarray}
\end{thm}
 \proof It is clear that
$$\ds\lim_{\alpha\rightarrow 0}r_{1}(\alpha)=-i(\gamma+(N/2)) .$$
In view of (\ref{v}), we deduce
$$\left|K\left(\frac{ix}{\cos \alpha},y\right)\right|\leq 1.$$
Then
$$\left|e^{\frac{i}{2}\tan (\alpha)(|x|^{2}+|y|^{2})}K\left(\frac{ix}{\cos \alpha},y\right)D_{k}f(y)\right|\leq |D_{k}f(y)|.$$
Let $y\in \R^{N}$ such that $|y|>1.$ Then
$$|D_{k}f(y)|\leq |y|^{2}|D_{k}f(y)|.$$
Since $y\longmapsto|y|^{2}D_{k}f \in L^{1}(\R^{N}, \ \omega_{k}(x)dx),$ it follows that $D_{k}f \in L^{1}(\R^{N}, \ \omega_{k}(x)dx)$ and the dominated convergence theorem implies
\begin{eqnarray*}
& \ & \ds\lim_{\alpha\rightarrow 0}r_{1}(\alpha) \ \frac{c_{k}}{2^{\gamma+(N/2)}}
\ds\int_{\R^{N}}e^{\frac{i}{2}\tan (\alpha)(|x|^{2}+|y|^{2})}K\left(\frac{ix}{\cos \alpha},y\right)D_{k}f(y)\omega_{k}(y)dy\\&=&-i(\gamma+(N/2))\frac{c_{k}}{2^{\gamma+(N/2)}}
\ds\int_{\R^{N}}K(ix,y)D_{k}f(y)\omega_{k}(y)dy\\&=&-i(\gamma+(N/2))D_{k}^{2}f(-x)\\
&=&-i(\gamma+(N/2))f(x) , \ a.e.
\end{eqnarray*}
From (\ref{dd}), we deduce
\begin{eqnarray*}
\ds\lim_{\alpha\rightarrow 0}r_{2}(\alpha,x,y)=\frac{i}{2}(|x|^{2}+|y|^{2})K(ix,y).
\end{eqnarray*}
By the previous Lemma, we have the following majorization:
\begin{eqnarray*}
\left|r_{2}(\alpha,x,y)D_{k}f(y)\right|\leq f_{1}(y)+f_{2}(y)+f_{3}(y),
\end{eqnarray*}
where
\begin{eqnarray*}
f_{1}(y)&=&\frac{1}{2}(1+\tan ^{2}\alpha_{0})|x|^{2}|D_{k}f(y)|,\\
f_{2}(y)&=&\frac{1}{2}(1+\tan ^{2}\alpha_{0})|y|^{2}|D_{k}f(y)|,\\
f_{3}(y)&=&\frac{|\sin (\alpha_{0})|}{\cos^{2} (\alpha_{0})} \ \sqrt{N} \ |x||y| |D_{k}f(y)|.
\end{eqnarray*}
Since $y\longmapsto|y|^{2}D_{k}f \in L^{1}(\R^{N}, \ \omega_{k}(x)dx),$ it follows that $f_{1}, f_{2}, f_{3} \in L^{1}(\R^{N}, \ \omega_{k}(x)dx)$ and therefore $f_{1}+f_{2}+f_{3} \in L^{1}(\R^{N}, \ \omega_{k}(x)dx).$
By virtue of the dominated convergence theorem, we have
\begin{eqnarray*}
& \  & \ds\lim_{\alpha\rightarrow 0}\frac{c_{k}}{2^{\gamma+\frac{N}{2}}}
\ds\int_{\R^{N}}r_{2}(\alpha,x,y)D_{k}f(y)\omega_{k}(y)dy\\&=& \frac{i}{2}\frac{c_{k}}{2^{\gamma+\frac{N}{2}}}\ds\int_{\R^{N}}(|x|^{2}+|y|^{2})K(ix,y)D_{k}f(y)\omega_{k}(y)dy\\
&=& \frac{i |x|^{2}}{2}\frac{c_{k}}{2^{\gamma+\frac{N}{2}}}\ds\int_{\R^{N}}K(ix,y)D_{k}f(y)\omega_{k}(y)dy\\&+&
\frac{i}{2}\frac{c_{k}}{2^{\gamma+\frac{N}{2}}}\ds\int_{\R^{N}}|y|^{2}K(ix,y)D_{k}f(y)\omega_{k}(y)dy\\&=&
\frac{i}{2}|x|^{2}f(x)+\frac{i}{2}\frac{c_{k}}{2^{\gamma+\frac{N}{2}}}\ds\int_{\R^{N}}|y|^{2}K(ix,y)D_{k}f(y)\omega_{k}(y)dy, a.e,\\&=&
\frac{i}{2}|x|^{2}f(x)+\frac{i}{2}D_{k}\left[|y|^{2}D_{k}(y)\right](-x) \ a.e.
\end{eqnarray*}
\begin{coro}\mbox{}\\
$1)$ $\mathcal{S}(\R^{N}) \subset \mathcal{W}\subset D(T).$\\
$2)$ For all $f\in \mathcal{S}(\R^{N}),$
\begin{eqnarray*}
-iTf=-(\gamma+(N/2))f+\frac{1}{2}(|x|^{2}-\Delta_{k})f
\end{eqnarray*}
\end{coro}
\proof \mbox{}\\
$1)$ Obvious.\\
$2)$ Let $f\in \mathcal{S}(\R^{N}).$ From Corollary 2.11 in \cite{Dunkl3}, we deduce
\begin{eqnarray*}
-y_{j}^{2}D_{k}f(y)=D_{k}[T_{j}^{2}f](y),
\end{eqnarray*}
where $j\in \left\{1,2,\dots,N\right\}.$
Then
\begin{eqnarray*}
-|y|^{2}D_{k}f(y)=D_{k}[\Delta_{k}f](y).
\end{eqnarray*}
Therefore
\begin{eqnarray}
-D_{k}\left[|y|^{2}D_{k}(y)\right](-x)&=&D_{k}^{2}[\Delta_{k}f(y)](-x)\nonumber\\&=&
\Delta_{k}f(x).\label{ee}
\end{eqnarray}
Finally, from (\ref{ff}) and (\ref{ee}) we obtain the desired result.
\begin{remark}It is clear that the operator $2iT-(2\gamma+N)$ is an extension on $\mathcal{W}$ of the Hermite operator $\mathcal{H}_{k}=\Delta_{k}-|x|^{2}$ studied by R\"{o}sler \cite{rosler} where it used another approach based on the notion of Lie algebra.\\
In the same context, we give a new proof of the following result established in \cite{rosler}
\end{remark}
\begin{coro} For $n\in \N$ and $p\in\mathcal{P}_{n},$ the function $f=e^{-\frac{|x|^{2}}{2}}e^{-\frac{\Delta_{k}}{4}}p(x)$ satisfies
\begin{eqnarray}
(\Delta_{k}-|x|^{2})f=-(2n+2\gamma+N)f.\label{z_{1}}
\end{eqnarray} In particular
\begin{eqnarray}
(\Delta_{k}-|x|^{2})h_{\nu}=-(2|\nu|+2\gamma+N)h_{\nu}.
\end{eqnarray}
\end{coro}
\proof Since $f\in \mathcal{S}(\R^{N}),$ (\ref{z_{1}}) is obtained by the use of the previous Corollary and (\ref{y_{1}})

\end{document}